\documentclass[a4paper,11pt]{amsart}
\usepackage{amssymb,amsfonts,amsxtra}
\usepackage[all]{xy}
\usepackage{fullpage}
\DeclareMathOperator{\hofib}{hofib}
\DeclareMathOperator{\Map}{Map}
\DeclareMathOperator*{\hocolim}{hocolim}

\def\lra{\longrightarrow}

\def\ra{\rightarrow}
\newtheorem{theorem}{Theorem}[section]
\newtheorem{cor}[theorem]{Corollary}
\newtheorem{lem}[theorem]{Lemma}
\newtheorem{prop}[theorem]{Proposition}

\theoremstyle{definition}
\newtheorem{defi}[theorem]{Definition}
\newtheorem{rem}[theorem]{Remark}
\numberwithin{equation}{section}
\begin{document}
\title[Cohomology theories for highly structured ring spectra]
{Cohomology theories for highly structured ring spectra}
\author{A. Lazarev}
\address{Mathematics Department, Bristol University, Bristol, BS8 1TW,
England.} \email{a.lazarev@bristol.ac.uk}
\begin{abstract}
This is a survey paper on cohomology theories for $A_\infty$ and
$E_\infty$ ring spectra. Different constructions and main
properties of topological Andr\'e-Quillen cohomology and of
topological derivations are described. We give sample calculations
of these cohomology theories and outline applications to the
existence of  $A_\infty$ and $E_\infty$ structures on various
spectra. We also explain the relationship between topological
derivations, spaces of multiplicative maps and moduli spaces of
multiplicative structures.
\end{abstract}
\maketitle
\section{Introduction} In recent years algebraic topology witnessed renewed
interest to highly structured ring spectra first introduced in
\cite{M1}. To a large extent it was caused by the discovery of a
strictly associative and symmetric smash product in the category
of spectra in \cite{EKMM}, \cite{HSS}. This allowed one to replace
the former highly technical notions of $A_\infty$ and $E_\infty$
ring spectra by equivalent but conceptually much simpler notions
of $S$-algebras and commutative $S$-algebras respectively.

An \emph{$S$-algebra} is just a monoid in the category of spectra
with strictly symmetric and associative smash product (hereafter
referred to as the category of \emph{$S$-modules}). Likewise, a
\emph{commutative $S$-algebra} is a commutative monoid in the
category of $S$-modules.

The most important formal property of categories of $S$-algebras
and commutative $S$-algebras is that both are \emph{topological
model categories} in the sense of Quillen, \cite{Q} as elaborated
in \cite{EKMM}. That means that together with the usual structure
of the closed model category (fibrations, cofibrations and weak
equivalences subject to a set of axioms) these categories are
topologically enriched, that is their $Hom$ sets are topological
spaces and the composition of morphisms is continuous. Moreover,
there exist \emph{tensors} and \emph{cotensors} of objects with
topological spaces that satisfy the usual adjunction isomorphisms
which hold for cartesian products and mapping spaces in the
category of topological spaces. In addition, this enrichment is
supposed to be compatible with the closed model structure in the
sense that an appropriate analogue of Quillen's corner axiom SM7
is satisfied.

This rich structure allows one to translate a lot of the notions
of conventional homotopy theory (homotopy relation, cellular
approximation, the formation of homotopy limits and colimits) in
the abstract setting. The language of closed model categories will
be freely used here and we refer the reader to the monograph
\cite{Hov} where the necessary details may be found.

The categories of $S$-algebras and commutative $S$-algebras are in
some sense analogous to the unstable category of topological
spaces. In these categories there are certain natural homotopy
invariant theories which play the role of singular cohomology for
topological spaces. In particular they are natural homes for
various obstruction groups. The corresponding theory is called
\emph{topological Andr\'e-Quillen cohomology} in the commutative
case and \emph{topological derivations} in the associative case.
We will use the abbreviation TAQ for the former and Der for the
latter. The purpose of the present paper is to give an overview of
recent results on these cohomology theories. The technical level
will be kept to a minimum, however we will try to outline proofs
of the results we formulate, especially when these proofs are easy
and conceptual. Occasionally the theorems we obtain are `new',
that is, have not appeared in print before; but in each case they
are either direct extensions of the known results or could be
obtained following similar patterns. Still, we tried to put these
results in different perspectives than those in the published
sources; in particular the associative and commutative cases are
treated as uniformly as possible.

The paper does not pretend to be comprehensive by any means. One
notable omission is the work of Robinson and Whitehouse on
$\Gamma$-homology of commutative rings \cite{RW} and the
subsequent work of Robinson \cite{Rob2} where $\Gamma$-homology
was used to show the existence of a commutative $S$-algebra
structure on $E_n$, the so-called Morava $E$-theory spectrum. One
should also mention the important papers \cite{Ric} and \cite{PR}
containing further results on $\Gamma$-homology. Our limited
expertise as well as the lack of space prevented us from including
these results in the survey. Another theme only briefly mentioned
here is the existence of the action of the Morava stabilizer group
on the spectrum $E_n$ (the Hopkins-Miller theorem). For this and
related topics we refer the reader to \cite{Rez}, \cite{L2} and
\cite{GH}.

The paper is organized as follows. In sections 2 and 3 we outline
the construction and basic properties of TAQ and Der. In section 4
we investigate the problem of lifting a map of $S$-algebras or
commutative $S$-algebras and introduce the notion of a
\emph{topological singular extension}. In section 5 we discuss an
alternative construction of TAQ and Der as \emph{stabilizations}
of appropriate forgetful functors. Section 6 presents calculations
of TAQ and Der of the Eilenberg-MacLane spectrum $H\mathbb{F}_p$.
In section 7 we show how the developed technology could be used to
produce $MU$-algebra structures on many complex oriented spectra.
Here $MU$ is the spectrum of complex cobordisms. In section 8 we
discuss the relationship between spaces of algebra maps and Der.
In section 9 we make a link between Der and spaces of
multiplicative structures on a given spectrum.

 \subsection*{Notation and conventions} The paper is written in the language
of $S$-modules of \cite{EKMM} and we adopt the terminology of the
cited reference. There is one exception: we use the terms
`cofibrant' and `cofibration' for what \cite{EKMM} calls
`$q$-cofibrant' and `$q$-cofibration'. This is because we never
have the chance to use cofibrations in the classical sense (maps
that satisfy the homotopy extension property). Except for Section
6 we work over an arbitrary (but fixed) cofibrant commutative
$S$-algebra $R$, smash products $\wedge$ and function objects
$F(-,-)$ mean $\wedge_R$ and $F_R(-,-)$; when we use a different
`ground' $S$-algebra $A$ this is indicated by a subscript such as
$\wedge_A$ and $F_A(-,-)$. The topological space of maps in a
topological category $\mathcal{C}$ is denoted by
$\Map_\mathcal{C}$. The category of unbased topological spaces is
denoted by ${\mathcal Top}$ and that of based spaces - by
${\mathcal Top}_*$. For an $S$-algebra $A$ the category of (left)
$A$-modules is denoted by $\mathcal{M}_A$. The category of
$A$-bimodules (which is the same as the category of $A\wedge
A^{op}$-modules where $A^{op}$ is the $R$-algebra $A$ with
opposite multiplication) will be denoted by $\mathcal{M}_{A-A}$.
The category of $R$-algebras is denoted by $\mathcal{C}^{ass}_R$
and that of \emph{commutative} $R$-algebras - by
$\mathcal{C}^{comm}_R$ The homotopy category of a closed model
category $\mathcal{C}$ is denoted by $h\mathcal{C}$ so, for
example $h\mathcal{C}^{comm}_R$ is the homotopy category of
commutative $R$-algebras. For an $R$-algebra $A$ we will denote by
$[-,-]_{A-mod}$ and $[-,-]_{A-bimod}$ the sets of homotopy classes
of $A$-module or $A$-bimodule maps respectively. The field
consisting of $p$ elements is denoted by $\mathbb{F}_p$.

 \subsection* {Acknowledgement} I would like to express my sincere gratitude
to Andy Baker and Birgit Richter for doing such a wonderful job of
organizing the Workshop on Structured Ring Spectra in Glasgow in
January 2002. I am also thankful to Maria Basterra for explaining
to me her joint work with Mike Mandell on TAQ and to Stefan
Schwede for pointing out that Theorem \ref{moduliderivation} is a
direct consequence of Theorem \ref{osn2}. Thanks are also due to
Mike Mandell, Paul Goerss, Nick Kuhn and Bill Dwyer for useful
discussions and comments made at various times.

\section{Topological Andr\'e-Quillen cohomology for commutative $S$-algebras}
Let $A$ be a commutative $R$-algebra which without loss of
generality will assumed to be cofibrant. Denote by
$\mathcal{C}^{comm}_{R/A}$ the category of commutative
$R$-algebras over $A$. An object of $\mathcal{C}^{comm}_{R/A}$ is
a commutative $R$-algebra $B$ supplied with an $R$-algebra map
$B\lra R$ (an augmentation). A morphism between two such objects
$B$ and $C$ is the following commutative diagram in
$\mathcal{C}^{comm}_R$:
\[\xymatrix{B\ar[rr]\ar[dr]&&C\ar[dl]\\&A&}\]
Then  $\mathcal{C}^{comm}_{R/A}$ inherits a topological model
category structure from $\mathcal{C}^{comm}_{R}$ so that a map in
$\mathcal{C}^{comm}_{R/A}$ is a cofibration  if it is so
considered as a map in $\mathcal{C}^{comm}_{R}$.  Note that in the
case $R=A$ the category
$\mathcal{C}^{comm}_{R/A}=\mathcal{C}^{comm}_{A/A}$ is pointed and
therefore is enriched over the category of pointed topological
spaces.

Let us denote by $\mathcal{N}_A$ the category of \emph{nonunital}
commutative $A$-algebras. An object in this category is an
$A$-module $M$ together with a strictly associative multiplication
map $M\wedge_A M\lra M$. The morphisms in $\mathcal{N}_A$ are
defined in the obvious fashion. Following \cite{Bas} we will refer
to an object of $\mathcal{N}_A$ as an $A$-NUCA. It is not hard to
prove that  $\mathcal{N}_A$ has a topological model structure
where weak equivalences are those maps which are weak equivalences
on underlying $A$-modules. The fibrations then are the maps that
are fibrations of underlying $A$-modules and the cofibrations are
the maps which have the LLP with respect to the acyclic fibration.
 Then we could form the homotopy category $h\mathcal{N}_{A}$ of
$\mathcal{N}_{A}$.

We want to show that the categories  $h\mathcal{C}^{comm}_{A/A}$
and $h\mathcal{N}_{A}$ are equivalent. In fact, more is true, see
Proposition \ref{qui}.

Let $K:\mathcal{N}_A\lra \mathcal{C}^{comm}_A$ denote the functor
which assigns to an $A$-NUCA $M$ the commutative $A$-algebra
$A\vee M$ with multiplication
\[(A\vee M)\wedge_A(A\vee M)\cong A\vee M\vee M\vee M \vee
M\wedge_A M\lra A\vee M\] given by the obvious maps on the first
three wedge summands and by the multiplication of $M$ on the last
one.

Clearly the $A$-algebra $A\vee M$ can be considered as an object
in $\mathcal{C}^{comm}_{A/A}$ via the canonical projection map
$A\vee M\lra A$. In other words the functor $K$ lands in fact in
the category $\mathcal{C}^{comm}_{A/A}$.

Now let $B$ be an $A$-algebra over $A$ and denote by $I$ its
`augmentation ideal', i.e. the fibre of the augmentation map
$B\lra A$. Then $I(B)$ is naturally an $A$-NUCA and we can
consider $I$ as a functor from $\mathcal{C}^{comm}_{A/A}$ to
$\mathcal{N}_{A}$.

Then we have the following \begin{prop}\label{qui} The functor $K$
is left adjoint to $I$. Moreover the functors $K$ and $I$
establish a Quillen equivalence between the categories
$\mathcal{N}_{A}$ and $\mathcal{C}^{comm}_{A/A}$.
\end{prop}
 To see that the
functors $K$ and $I$ are adjoint notice that the category
$\mathcal{N}_{A}$ is in fact a category of algebras over a certain
monad $\mathbb{A}$ in $A$-modules. This monad is specified by
$\mathbb{A}M:=\bigvee_{i>0}M^{\wedge_A i}/\Sigma_i$ where
$\Sigma_i$ is the symmetric group on $i$ symbols. Further note
that $K(\mathbb{A}M)$ is a free commutative $A$-algebra on $M$.
Therefore if $M=\mathbb{A}M$ for an $A$-module $M$ we have
\[\mathcal{C}^{comm}_{A/A}(K(\mathbb{A}M),B)\cong
\mathcal{M}_A(M,I(B))\cong \mathcal{N}_A(\mathbb{A}M,B).\] To get
the adjointness isomorphism for a general $A$-NUCA $X$ it suffices
to notice that there is a canonical split coequalizer (the
beginning of the monadic bar-construction for $M$):
\[\xymatrix{\mathbb{A}\mathbb{A}X\ar@<-0.5ex>[r]\ar@<0.5ex>[r]&\mathbb{A}X\ar[r]&M}.\]
Explicitly, if $g:M\lra I(B)$ is a morphism in $\mathcal{N}_A$
then its adjoint $\tilde{g}:{A\vee M\lra B}$ is the composite
\[\xymatrix{\tilde{g}:A\vee M\ar[r]^{id\vee g}&A\vee I(B)\ar[r]^-{1\vee i}&B}\]
where $i$ is the canonical map $I(B)\lra B$.

 Further straightforward arguments show that the functor $K$
preserves cofibrations and acyclic cofibrations and that the
adjunction described above determines an equivalence on the level
of homotopy categories.

We will denote the homotopy invariant extension (also called a
\emph{total derived functor}) of the functor $I$ by $\mathbf{R}I$.
For an object $B$ in $\mathcal{C}^{comm}_{A/A}$ denote by
$\tilde{B}$ its fibrant replacement. Then
$\mathbf{R}I(B)=I(\tilde{B})$.
 Next we define the functor of `taking the indecomposables'
that assigns to an $A$-NUCA $N$ the $A$-module $Q(N)$ given by the
cofibre sequence:
\[\xymatrix{N\wedge_AN\ar[r]^-m&N\ar[r]&Q(N)}.\]
Here $m$ stands for the multiplication map. The functor $Q$ has a
right adjoint functor $Z:\mathcal{M}_A\lra \mathcal{N}_A$ which is
given by considering $A$-modules as objects in $\mathcal{N}_A$
with zero multiplication. It is easy to see that this adjunction
passes to the homotopy categories. Let $\mathbf{L}Q$ denote the
total derived functor of $Q$. Explicitly,
$\mathbf{L}Q(N)=Q(\tilde{N})$ where $\tilde{N}$ is the cofibrant
replacement of the $A$-NUCA $N$.

It is clear that $Q(\mathbb{A}N)\cong N$. However for a general
cofibrant $A$-NUCA the functor $Q(N)$ is very hard to compute. One
can approach its computation as follows. The functor
$\mathbb{A}:\mathcal{M_A}\lra \mathcal{M}_A[\mathbb{A}]$ from
$A$-modules to $\mathbb{A}$-algebras is left adjoint to the
forgetful functor $U:\mathcal{M}_A[\mathbb{A}]\lra \mathcal{M}_A$.
There results a monad in $\mathcal{M}_A$. Given an $A$-NUCA $N$
denote by $B_*(N)=B_*(\mathbb{A},U\mathbb{A},UN)$ its monadic
bar-construction. Explicitly $B_*N$ is the cosimplicial $A$-NUCA
with $B_n(N)=\mathbb{A}^{n+1}N$ the faces and codegeneracies are
the standard ones, cf. \cite{May}. The geometric realization
$|B_*(N)|$ of $B_*(N)$ is weakly equivalent to $N$. Then one has
the following
\begin{prop}\label{har} For a cofibrant $A$-NUCA $N$ one has the following
weak equivalence of $A$-modules: $Q(|B_*(N)|)\simeq Q(N)$.
\end{prop}
Despite the innocent appearance of this proposition it is not at
all obvious. The problem is that the functor $Q$ only preserves
weak equivalences between  cofibrant objects and $|B_*(N)|$ is not
a cofibrant $A$-NUCA. However in the end it turns out that it is
close enough to being cofibrant to give the result; much of the
work \cite{Bas} is devoted to overcoming this point.

Further associated to the simplicial $A$-NUCA $B_*(N)$ is the
spectral sequence which computes $\pi_*|B_*(N)|=\pi_*N$:
\begin{equation} \label{sp1}E^1_{i,j}=\pi_iQ(\mathbb{A}^{j+1}N)=
\pi_i\mathbb{A}^jN\Longrightarrow \pi_{i-j}Q(N).\end{equation}

 We now have all the
ingredients to define the \emph{abelianization} functor:
\[\Omega^{comm}_A:h\mathcal{C}^{comm}_{R/A}\lra h\mathcal{M}_A.\] \begin{defi}Let
$C$ be a commutative $R$-algebra over $A$ and $M$ be a $A$-module.
Then $\Omega_A^{comm}(C):=\mathbf{L}Q\circ
\mathbf{R}I(C\wedge_R^{\mathbf{L}} A)$.
\end{defi}
(Note that $C\wedge_R^{\mathbf{L}} A$ is an object of
$\mathcal{C}^{comm}_{A/A}$ in an obvious way.) Of course we only
defined the abelianization functor on objects, but its extension
to morphisms is immediate. The main property of the abelianization
functor is that it is left adjoint to the `square-zero' extension
functor: if $M$ is a $A$-module then the square zero extension of
$A$ by $M$ is the $A$-algebra over $A$ whose underlying $R$-module
is $A\vee M$ with the obvious multiplication. More precisely we
have the following
\begin{theorem}\label{osn1} There is the following natural isomorphism:
\[h\mathcal{C}^{comm}_{R/A}(C,A\vee M)\cong h\mathcal{M}_A (\Omega^{comm}_A(C),M)\]
where $C$ is a commutative $R$-algebra over $A$ and $M$ is a
$A$-module.\end{theorem} The proof of the theorem is rather formal
and relies on the already established properties of the functors
$I,K,Z$ and $Q$.

\begin{rem}We would like to emphasize here that the isomorphism of
Theorem \ref{osn1} holds on the level of \emph{homotopy}
categories and is not a reflection of an adjunction between strict
categories. Indeed, even though the abelianization functor could
be considered as a point-set level functor its definition involves
composition of right adjoint and left adjoint functors. Therefore
one cannot expect any good formal properties of $\Omega^{comm}_A$
on the point-set level. The situation improves upon passing to
homotopy categories because one of these functors, namely $RI$,
becomes an equivalence.
\end{rem}
 We will sometimes need the enriched version of Theorem
\ref{osn1}. The category $\mathcal{C}^{comm}_{R/A}$ is enriched
over unbased topological spaces. For two objects $B,A$ of a
topological category $\mathcal C$ the space of maps
$\Map_{\mathcal{C}/A}(B,C)$ from $B$ to $C$ \emph{over} $A$ is
defined from the pullback diagram
\[\xymatrix{\Map_{\mathcal{C}/A}(B,C)\ar[r]\ar[d]&pt\ar[d]\\
 \Map_{\mathcal{C}}(B,C)\ar[r]&
\Map_{\mathcal{C}}(B,A)}\] where the right downward arrow is just
picking the structure map $B\lra A$ in $\Map_{\mathcal{C}}(B,A)$.
It is easy to see that the cotensor $(A\vee M)^X$ of $X$ and
$A\vee M$ in $\mathcal{C}^{comm}_{R/A}$ is weakly equivalent to
$A\vee M^X$. Then for any $CW$-complex $X$ we have the following
isomorphisms:
\begin{align*}h{\mathcal Top}(X,\Map_{\mathcal{M}_C}(\Omega^{comm}_A(C),M))&\cong
h\mathcal{M}_C(\Omega^{comm}_A(C),M^X)\\ &\cong
h\mathcal{C}^{comm}_{R/A}(C, A\vee M^X)\\ &\cong h{\mathcal
Top}(X,\Map_{\mathcal C_{R/A}}(C,A\vee M)).
\end{align*}
Therefore the topological spaces $\Map_{\mathcal C_{R/A}}(C,A\vee
M)$ and $\Map_{\mathcal{M}_C}(\Omega^{comm}_A(C),M)$ are weakly
equivalent.

We will denote the $A$-module $\Omega^{comm}_A(A)=\mathbf{L}Q\circ
\mathbf{R}I(A\wedge^{\mathbf{L}}_R A)$ simply by
$\Omega^{comm}_{A}$ so that
\[h\mathcal{C}^{comm}_{R/A}(A,A\vee M)\cong
h\mathcal{M}_A(\Omega^{comm}_{A},M).\]  The set
$h\mathcal{C}^{comm}_{R/A}(A,A\vee M)$ could be interpreted as
derivations of the commutative $S$-algebra $A$ with coefficients
in the $A$-module $M$. In particular we see, that the set of such
derivations is an abelian group (note that apriori the set of
homotopy classes of (commutative) $S$-algebra maps does not carry
a structure of a group, let alone an abelian group).

We can now define the \emph{topological Andr\'e-Quillen
cohomology} spectrum of $A$ relative to $R$ with coefficients in
$M$ by
\[{\bf TAQ}_R(A,M):=F_A(\Omega^{comm}_{A},M)\]
and similarly the \emph{topological Andr\'e-Quillen homology}
spectrum of $A$ relative to $R$ with coefficients in $M$:
\[{\bf TAQ}^R(A,M):=\Omega^{comm}_{A}\wedge_AM.\]
We will refer to the homotopy groups of these spectra as
topological Andr\'e-Quillen (co)homology of $A$ relative to $R$:
\[{TAQ}_R^*(A,M):=\pi_{-*}F_A(\Omega^{comm}_{A},M);\]
\[TAQ^R_*(A,M):=\pi_*\Omega^{comm}_{A}\wedge _AM.\]
\begin{rem}\label{inf1}The enriched version of Theorem \ref{osn1} gives us the  weak
equivalence of spaces $\Map_{\mathcal{C}^{comm}_{R/A}}(B,A\vee
M)\simeq \Map_{\mathcal{M}_A}(\Omega^{comm}_A(C),M)$. The latter
mapping space is in turn weakly equivalent to the zeroth space of
the spectrum $F_A(\Omega^{comm}_A,M)={\bf TAQ}_R(A,M)$. We obtain
the weak equivalence
\[\Map_{\mathcal{C}^{comm}_{R/A}}(B,A\vee M)\simeq \Omega^\infty
{\bf TAQ}_R(A,M).\]\end{rem}

We are most interested in topological Andr\'e -Quillen
\emph{cohomology} since they are related to the obstruction
theory. Repeat, that for a commutative $R$-algebra $A$ and an
$A$-module $M$ the abelian group $TAQ_R^0(A,M)$ is identified with
the set of \emph{commutative derivations} of $A$ with values in
$M$ that is the homotopy classes of maps in
$\mathcal{C}^{comm}_{R/A}$ from $A$ into $A\vee M$. We say
`commutative derivations' to distinguish them from just
`derivations' or `topological derivations' which will be
considered in the next section and refer to the \emph{bimodule}
derivations in the context of associative $R$-algebras.

Let us now discuss the fuctoriality of TAQ. Since by definition
${\bf TAQ}_R(A,M)=F_{A}(\Omega^{comm}_A,M)$  we see that it is
covariant with respect to the variable $M$. On the level of
commutative derivations it can be seen as follows. A map of
$A$-bimodules $f:M\lra N$ determines a map in
$\mathcal{C}^{comm}_{R/A}:A\vee M\lra A\vee N$. Then a topological
derivation $d:A\lra A\vee M$ determines a topological derivation
$f_*d:A\lra A\vee M\lra A\vee N$.

If $g:X\lra A$ is a map of $R$-algebras then the composite map
\[\xymatrix{X\ar^g[r]& A\ar^-d[r] &A\vee M}\] is a map in
$\mathcal{C}^{comm}_{A/R}$. Therefore by Theorem \ref{osn1} it
corresponds to a map of in
$h\mathcal{M}_{A}:\Omega^{comm}_A(X)\lra M$ which is the same as a
map in $h\mathcal{M}_{X}:\Omega^{comm}_X\lra M$. The latter
corresponds, again by Theorem \ref{osn1} to a commutative
derivation of $X$ with values in the $X$-module $M$. We will
denote this derivation by  $g^*d:X\lra X\vee M$.

We can relate the topological Andr\'e-Quillen cohomology to
ordinary $R$-module cohomology as follows. For a cofibrant
$R$-algebra $A$ consider the composite map of $R$-modules
\[l_R:A\lra A\vee\Omega^{comm}_{A}\lra\Omega^{comm}_{A}.\]
Here the first map is the $R$-algebra map adjoint to the identity
on $\Omega^{comm}_{A}$ and the second map is the projection onto
the wedge summand. The map $l_R$ induces a forgetful map from TAQ
to $R$-module cohomology:
\[l^*_R:TAQ_R^*(A,M)=[\Omega^{comm}_{A},M]^*_{A-mod}\lra
[\Omega^{comm}_{A},M]^*_{R-mod}\lra [A,M]^*_{R-mod}.\] To describe
the image of $l^*_R$ let us introduce the notion of a
\emph{primitive operation} $E\lra M$ for a $R$-ring spectrum $E$
and an $E$-bimodule spectrum $M$. Namely, the map $p\in
[E,M]^*_{R-mod}$ is primitive if the following diagram is
commutative in $h\mathcal{M}_R$:
\[\xymatrix{E\wedge E\ar[r]^m\ar[d]_{1\wedge p\vee p\wedge
1}&E\ar[d]^p\\ E\wedge M\vee M\wedge E\ar[r]^-{m_l\vee m_r}&M}\]
where $m_l$ and $m_r$ denote the left and right actions of $E$ on
$M$. (In other words, $p$ is a derivation up to homotopy but we
refrain from using this term to avoid confusion.) Then it is easy
to see that the image of $l^*_R$ is contained in the subspace of
primitive operations in $[A,M]^*_{R-mod}$. Of course the left and
right $A$-module spectrum structures on $M$ coincide in the
commutative case.

\section{Topological derivations}
In this section we construct the analogue of TAQ for not
necessarily commutative $S$-algebras. Traditionally the
Andr\'e-Quillen cohomology was considered for commutative algebras
only. Therefore we will reserve this term for the commutative case
and call its analogue for associative $S$-algebras
\emph{topological derivations}. The main result we are going to
describe here is the analogue of Theorem \ref{osn1} in the context
of noncommutative $S$-algebras. Note, however, that the
construction and the proof differ considerably from the
commutative case.

Let $A$ be a cofibrant $R$-algebra. We define the module of
differentials $\Omega^{ass}_{A}$ from the following homotopy fibre
sequence
\begin{equation}\label{dif}\xymatrix{\Omega^{ass}_{A}\ar[r]&A\wedge
A\ar[r]&A}.
\end{equation}
Here the second arrow is the multiplication map. Clearly
$\Omega^{ass}_{A}$ is an $A$-bimodule.  Note that the sequence
(\ref{dif}) splits in the homotopy category of left $A$-modules by
the map \[\xymatrix{A\ar^-\cong[r]& A\wedge R\ar^-{id\wedge
1}[r]&A\wedge A}\] and therefore the $A$-bimodule $\Omega^{ass}_A$
is equivalent as a left $A$-module to $A\wedge A/R$. Similarly
$\Omega^{ass}_A$ is equivalent as a \emph{right} $A$-module to
$A/R\wedge A$.

Denote the category of objects over $A$ inside
$\mathcal{C}^{ass}_{R}$ by $\mathcal{C}^{ass}_{R/A}$. The category
$\mathcal{C}^{ass}_{R/A}$ has a topological model category
structure inherited from $\mathcal{C}^{ass}_{R}$.

Let $M$ be an $A$-bimodule. Then we could form the $R$-algebra
with the underlying $R$-module $A\vee M$, the `square-zero
extension' of $A$ by $M$. There is an obvious product and
augmentation maps making $A\vee M$ into an $R$-algebra over $A$.

We now introduce the analogue of the abelianization functor in the
noncommutative context. Let $B$ be a cofibrant $R$-algebra over
$A$ and denote by $\Omega^{ass}_A(B)$ the $A$-bimodule
\[A\wedge_B\Omega^{ass}_B\wedge_BA\cong A\wedge A^{op}\wedge_{B\wedge
B^{op}}\Omega^{ass}_B.\] Notice that $\Omega^{ass}_A(B)$ is an
$A$-bimodule.
\begin{theorem}\label{osn} There is a natural equivalence
\[h\mathcal{C}^{ass}_{R/A}(B, A\vee M)\cong
h\mathcal{M}_{A-A}(\Omega^{ass}_A(B),M).\] \end{theorem} Let us
show how to associate to any map $B\lra A\vee M$ in
$ho\mathcal{C}^{ass}_{R/A}$ a map of $A$-bimodules
$\Omega^{ass}_A(B)\lra M$.

Denote by $\overline{A\vee M}$ the fibrant cofibrant approximation
of $A\vee M$ in the category of $R$-algebras over $A$. It suffices
to construct a `universal' map of $\overline{A\vee M}$-bimodules
\begin{equation}\label{oh}\Omega^{ass}_A(\overline{A\vee M})\lra M.
\end{equation} Indeed the for a map $B\lra \overline{A\vee M}$ in
$\mathcal{C}^{ass}_{R/A}$ we have a composite map
\[\Omega^{ass}_A(B)\lra \Omega^{ass}_A(\overline{A\vee M})\lra M\] and
therefore a correspondence
\begin{equation}\label{ohh}h\mathcal{C}^{ass}_{R/A}(B,\overline{A\vee
M})=h\mathcal{C}^{ass}_{R/A}(B,A\vee M)\lra
h\mathcal{M}_{A-A}(\Omega^{ass}_A(B),M)\end{equation} as desired.

Further since \[\Omega^{ass}_A(\overline{A\vee
M})=A\wedge_{\overline{A\vee M}}\Omega^{ass}_{\overline{A\vee
M}}\wedge_{\overline{A\vee M}}A\]  a map of (\ref{oh}) is is the
same as a map of $\overline{A\vee M}$-bimodules
$\Omega^{ass}_{\overline{A\vee M}}\lra M$. Instead of the latter
map we construct a map $\Omega^{ass}_{A\vee M}\lra M$. This would
be good enough since even $A\vee M$ is not a cofibrant object in
$\mathcal{C}^{ass}_{R/A}$ the smash product $(A\vee M)\wedge
(A\vee M)$ clearly represents the derived smash product and
therefore $\Omega^{ass}_{\overline{A\vee M}}$ is weakly equivalent
to $\Omega^{ass}_{A\vee M}$.
 Consider the following diagram in the homotopy category of $A\vee
 M$-bimodules:
 \[\xymatrix{A\wedge A\vee A\wedge M\vee M\wedge A\vee M\wedge
 M\ar[r]\ar[d]&A\vee M\ar@{=}[d]\\A\vee M\vee M\ar[r]&A\vee M}\]
Here the left vertical arrow is determined by the $R$-algebra
structure on $A$, an $A$-bimodule structure on $M$ and is zero on
the last summand. The lower horizontal arrow is zero on $A$,
identity on the first $M$-summand and minus identity on the last
$M$-summand. Then the homotopy fibre of the upper row is
equivalent to $\Omega^{ass}_{A\vee M}$ by definition and the
homotopy fibre of the lower row is equivalent to $M$. There
results a map of $A\vee M$-bimodules $\Omega^{ass}_{A\vee M}\lra
M$ as desired.

The proof of the theorem is then completed by showing that the
above correspondence (\ref{ohh}) is in fact an equivalence if $B$
is equal to the free algebra \[B=TV=R\vee V\vee V^{\otimes
2}\vee\ldots\] where $V$ is a cofibrant $R$-module over $A$ and
then resolving a general $B$ by a monadic bar-construction. Note
that the homotopy fibre of the multiplication map $TV\wedge TV\lra
TV$ is equivalent to the $TV$-bimodule $TV\wedge V\wedge TV$ and
therefore
\[\Omega^{ass}_{TV}\simeq TV\wedge V\wedge TV.\]

We have the `universal derivation' map $A\lra A\vee
\Omega^{ass}_A$ which is adjoint to the identity map
$\Omega^{ass}_A(A)=\Omega^{ass}_A\lra \Omega^{ass}_A$. Recall that
as an $R$-module $\Omega^{ass}_A$ (even as a left $A$-module) is
weakly equivalent to $A\wedge A/R$. Taking the composition of the
universal derivation with the canonical projection $A\vee
\Omega^{ass}_A\lra \Omega_A$ we get a map $\phi:A\lra
\Omega^{ass}_A\simeq A\wedge A/R$. Then as a map of $R$-modules
$\phi$ is homotopic to the map \[\xymatrix{A\cong R\wedge
A\ar^{1\wedge d}[r]&A\wedge A/R}\] where $d$ is is defined as the
second arrow in the homotopy cofibre sequence
\[\xymatrix{R\ar^1[r]&A\ar^d[r]&A/R.}\] This could be seen by
first taking $A$ to be the free algebra $TV$ and then resolving a
general $R$-algebra $A$ by a simplicial construction consisting of
free $R$-algebras.

Let us denote by $Der^0_R(A,M)$ the set of homotopy classes of
maps $A\lra A\vee M$ in $h\mathcal{C}^{ass}_{A/R}$. We will call
elements of $Der^0_R(A,M)$ \emph{topological derivations} of $A$
with values in the $A$-bimodule $M$. Then Theorem \ref{osn} shows
that $Der^0_R(A,M)$ is in fact an abelian group. Moreover,
$Der^0_R(A,M)$ is the zeroth homotopy group of the function
spectrum $F_{A\wedge A^{op}}(\Omega^{ass}_A,M)$. We will call it
the spectrum of topological derivations of $A$ with values in $M$
and denote it by ${\bf Der}_R(A,M)$. The homotopy groups of ${\bf
Der}_R(A,M)$ will be denoted by ${Der}^i_R(A,M)$, so
${Der}^i_R(A,M)=\pi_{-i}{Der}^i_R(A,M)$.

\begin{rem}\label{inf} Arguing in the same way as for commutative derivations
we see that there is a weak equivalence of topological spaces
\[\Map_{\mathcal{C}^{ass}_{R/A}}(A,A\vee M)\simeq \Omega^\infty{\bf
Der}_R(A,M).\] In particular the space
$\Map_{\mathcal{C}^{ass}_{R/A}}(A,A\vee M)$ is an infinite loop
space.
\end{rem}

Let $\tilde{A}$ be a cofibrant replacement of $A$ as an
$A$-bimodule. Replacing $A$ with $\tilde{A}$ in the homotopy fibre
sequence (\ref{dif}) and applying the functor $F_{A\wedge
A^{op}}(?,M)$ to it we obtain the homotopy fibre sequence:
\[F_{A\wedge A^{op}}(\tilde{A},M)\lra M\lra{\bf Der}_R(A,M).\] The function
spectrum $F_{A\wedge A^{op}}(\tilde{A},M)$ is called
\emph{topological Hochschild cohomology} spectrum of $A$ with
values in $M$ and has a special notation ${\bf THH}_R(A,M)$. So we
obtain the following homotopy fibre sequence relating topological
Hochschild cohomology and topological derivations:
\begin{equation}\label{hoch}{\bf THH}_R(A,M)\lra M\lra{\bf
Der}_R(A,M).\end{equation} This fibre sequence to a large extent
reduces the calculation of ${Der}^*_R(A,M)$ to that of topological
Hochschild cohomology groups ${THH}^*_R(A,M):=\pi_{-i}{\bf
THH}_R(A,M)$. The latter groups are relatively computable due to
the hypercohomology spectral sequence of \cite{EKMM}:
\[Ext^{**}_{(A\wedge A^{op})_*}(A_*,M_*)\Longrightarrow {
THH}^*_R(A,M).\] Similarly to the commutative case there is a
forgetful map
\[\xymatrix{l^*_R:{Der}^*_R(A,M)\ar[r]&[A,M]^*_{R-mod}}\]
whose image is contained in the subspace of primitive operations
in $[A,M]^*_{R-mod}$.
\begin{rem}\label{dual} There is also the notion of topological Hochschild
\emph{homology} spectrum \[{\bf
THH}^R(A,M):=\tilde{A}\wedge_{A\wedge A^{op}}M.\] Topological
Hochschild homology plays an important role in computations of
algebraic $K$-theory of rings, cf. \cite{Madsen}. If the
$R$-algebra $A$ is commutative then both ${\bf THH}_R$ and ${\bf
THH}^R$ are $R$-modules and the is a duality between them in
$h\mathcal{M}_A$:
\[{\bf THH}_R(A,M)\simeq F_A({\bf THH}^R(A,M),A).\]\end{rem}

\begin{rem}Just as topological Andr\'e-Quillen cohomology
the spectrum ${\bf Der}_R(A,M)$ is covariant with respect to $M$
and contravariant with respect to $A$. The arguments are precisely
the same as in the commutative context except that we refer to
Theorem \ref{osn} instead of Theorem \ref{osn1} and $A$-modules
are replaced with $A$-bimodules. For a map of $A$-bimodules
$f:M\lra N$ and a topological derivation $d:A\lra A\vee M$ we have
a topological derivation $f_*d$ of $A$ with values in $N$. For a
map of $R$-algebras $g:X\lra A$ we obtain a topological derivation
$g^*d$ of $X$ with values in $M$.\end{rem} Sometimes in the
noncommutative context it is useful to consider the so-called
\emph{relative topological derivations} which we will now define.
Let $R^\prime$ be a not necessarily commutative $R$-algebra.
Denote by $\mathcal{C}^{ass}_{R^{\prime}}$ the category of
$R$-algebras $A$ supplied with a fixed $R$-algebra map
$R^\prime\lra A$, not necessalily central. In other words
$\mathcal{C}^{ass}_{R^{\prime}}$ is the undercategory of
$R^\prime$ in $\mathcal{C}_{R}$. The objects of
$\mathcal{C}^{ass}_{R^{\prime}}$ will be called
\emph{$R^\prime$-algebras} by a slight abuse of language. Let
$m:A\wedge A\lra A$ be the multiplication map. Note that there
exists a unique map $m^\prime:A\wedge_{R^\prime} A\lra A$ such
that the diagram
\[\xymatrix{A\wedge
A\ar^m[rr]\ar[dr]&&A\\&A\wedge_{R^\prime}A\ar^{m^\prime}[ur]}\] is
commutative.

Now suppose without loss of generality that the $R$-algebra
$R^\prime $ is cofibrant and that the structure map $R^\prime\lra
A$ is a cofibration of $R$-algebras. We define  the
$R^\prime$-relative module of differentials
$\Omega^{ass}_{A/R^\prime}$ from the following homotopy fibre
sequence:
\[\xymatrix{\Omega^{ass}_{A/R^\prime}\ar[r]&A\wedge_{R^\prime} A\ar^-{m^\prime}[r]& A}.\]
Let $\mathcal{C}^{ass}_{A/R^{\prime}}$ denote the category whose
objects are $R^\prime$-algebras $B$ supplied with an
$R^\prime$-algebra map $B\lra A$. Note that if $M$ is an
$A$-bimodule then $A\vee M$ is an object of
$\mathcal{C}^{ass}_{A/R^{\prime}}$. This gives a functor
$\mathcal{M}_{A-A}\lra \mathcal{C}^{ass}_{A/R^{\prime}}$. It turns
out that on the level of homotopy categories this functor admits a
left adjoint \[B\lra
\Omega^{ass}_{A/R^{\prime}}(B):=A\wedge_B\Omega^{ass}_{B/R^{\prime}}\wedge_BA.\]
More precisely there is a natural equivalence
\[h\mathcal{C}^{ass}_{A/R^{\prime}}(B,A\vee M)\cong
h\mathcal{M}_{A-A}(\Omega^{ass}_{A/R^{\prime}}(B),M).\] The proof
is similar to the one in the absolute case. One makes use of the
`free' $R^\prime$-algebra  functor $T_{R^\prime}$ given for an
$R^\prime$-bimodule $V$ as \[T_{R^\prime}V=R^\prime\vee (V\vee
V\wedge_{R^\prime}V)\vee\ldots\] instead of the usual free
$R$-algebra on an $R$-module $V$.

The homotopy classes of maps in
$h\mathcal{C}^{ass}_{A/R^\prime}(A,A\vee M)$ are called
topological derivations of $A$ with values in $M$ \emph{relative
to $R^\prime$}. An algebraic analogue of relative topological
derivations are the derivations of an algebra $A$ with
coefficients in an $A$-bimodule $M$ {vanishing on a subalgebra
$R^\prime$ of $A$}. Relative modules of differentials were studied
in the context of noncommutative differential forms in \cite{QC}.
Replacing $\Omega^{ass}_{R/R^{\prime}}$ with
$\tilde{\Omega}^{ass}_{R/R^{\prime}}$, its cofibrant approximation
as an $A$-bimodule we define the $R$-module of
$R^{\prime}$-relative derivations of $A$ with values in an
$A$-bimodule $M$:\[{\bf Der}_{R/R^{\prime}}(A,M):=F_{A\wedge
A^{op}}(\tilde{\Omega}^{ass}_{R/R^{\prime}},M).\] It is easy to
see that there is the following homotopy fibre sequence of
$R$-modules:
\[\xymatrix{{\bf THH}_R(A,M)\ar[r]&{\bf THH}_R(R^\prime,M)\ar[r]&{\bf
Der}_{R/R^{\prime}}(A,M)}\] which specializes to (\ref{hoch}) when
$R=R^\prime$ .
\section{Obstruction theory}
In this section we show that topological Andre-Quillen cohomology
in the commutative case and topological derivations in the
associative case are a natural home for lifting algebra maps. For
definiteness we work in the context of associative algebras,
however the results have obvious analogues, with the same proofs,
in the commutative case. Topological derivations are replaced with
commutative derivations, spectra ${\bf Der}_R(-,-)$ - with spectra
${\bf TAQ}_R(-,-)$ and the spaces
$\Map_{\mathcal{C}^{ass}_R}(-,-)$ - with spaces
$\Map_{\mathcal{C}^{comm}_R}(-,-)$.

 Let $A$ be a cofibrant $R$-algebra, $I$
be an $A$-bimodule and $d:A\lra A\vee \Sigma I$ be a derivation of
$A$ with coefficients in $\Sigma I$, the suspension of $I$. Let
$\epsilon:A\lra A\vee \Sigma A$ be the inclusion of a retract.
Define the $R$-algebra $B$ from the following homotopy pullback
square of $R$-algebras:
\[\xymatrix{B\ar[d]\ar[r]&A\ar^{\epsilon}[d]\\A\ar^-d[r]&A\vee \Sigma
I}\] It is clear that the homotopy fibre of the map $A$ is weakly
equivalent to the $A$-bimodule $I$. \begin{defi} The homotopy
fibre sequence \[I\lra B\lra A\] is called the topological
singular extension of $A$ by $I$ associated with the derivation
$d:A\lra A\vee \Sigma I$.\end{defi}
\begin{theorem}\label{obs} Let $I\rightarrow B\rightarrow A$ be a singular extension of $R$-algebras
associated with a derivation $d:B\rightarrow B\vee\Sigma I$ and
$f:X\rightarrow A$ a map of $R$-algebras where the $R$-algebra $X$
is cofibrant. Then: \begin{enumerate} \item The map $f$ lifts to
an $R$-algebra map $X\rightarrow B$ if and only if the induced
derivation $f^*d\in Der^1_R(X,I)=Der^0_R(X,\Sigma I)$ is homotopic
to zero.\item Assuming that a lifting exists in the fibration
\begin{equation}\label{fib}\Map_{\mathcal{C}^{ass}_R}(X,B)\rightarrow
\Map_{\mathcal{C}^{ass}_R}(X,A)
\end{equation}
the homotopy fibre over the point $f\in
\Map_{\mathcal{C}^{ass}_R}(X,A)$ is weakly equivalent to
$\Omega^\infty {\bf Der}_R(X,I)$ (the $0$th space of the spectrum
${\bf Der}_R(X,I))$.\end{enumerate}
\end{theorem}
The proof is rather formal. Assume first that a lifting of $f$
exists. Then
  $d\circ f$ factors as
\[\xymatrix{X\ar^f[r]& B\ar[r]&A\ar^\epsilon[r]&A\vee
\Sigma I}\] which means that the derivation $d\circ f$ is trivial,
Conversely, if $d\circ f$ is trivial then the diagram of
$R$-algebras \[\xymatrix{ X\ar^f[r]\ar^f[d]&A\ar^\epsilon[d]\\
A\ar^d[r]&A\vee \Sigma I }\] commutes up to homotopy and by the
universal property of the homotopy pullback there is a map
$\tilde{f}:X\rightarrow B$ lifting $f$.

 For the second part suppose that the lifting of $f$ exists, so the (homotopy) fibre
of the map (\ref{fib}) is nonempty. We have the following diagram
of $R$-algebras. \[\xymatrix{
B\ar[dd]\ar[rr]\ar[dr]&&A\ar[dd]\ar@{=}[dl]\\&A&\\A\ar@{=}[ur]\ar^d[rr]&&
A\vee \Sigma I\ar[ul] }\]
  Changing $B$ in its homotopy class if necessary we can
arrange that this diagram be strictly commutative.
 Notice
that the outer square is a homotopy pullback of $R$-algebras.
Applying the functor $\Map_{\mathcal{C}^{ass}_R}(X,-)$ to this
diagram we get the diagram of spaces \[\xymatrix{
\Map_{\mathcal{C}^{ass}_R}(X,B)\ar[dd]\ar[rr]\ar[dr]&&\Map_{\mathcal{C}^{ass}_R}(X,A)
\ar[dd]\ar@{=}[dl]\\&\Map_{\mathcal{C}^{ass}_R}(X,A)&\\
\Map_{\mathcal{C}^{ass}_R}(X,A)\ar@{=}[ur]\ar[rr]&&
\Map_{\mathcal{C}^{ass}_R}(X,A\vee \Sigma I)\ar[ul] }\] Again the
outer square is a homotopy pullback (of spaces). Taking the
homotopy fibres of the maps from the outer square to the center
(over $f\in \Map_{\mathcal{C}_{R}^{ass}}$) we get the following
homotopy pullback of spaces:
\begin{equation}\label{usd}\xymatrix{
\hofib \Map_{\mathcal{C}^{ass}_R}(X,B)\ra
\Map_{\mathcal{C}^{ass}_R}(X,A)\ar[r]\ar[d]&pt\ar[d]\\
 pt\ar[r]&\hofib \Map_{\mathcal{C}^{ass}_R}(X,A\vee \Sigma I)
 \rightarrow \Map_{\mathcal{C}^{ass}_R}(X,A) }
\end{equation}
Notice that the space in the right lower corner of (\ref{usd}) is
canonically weakly equivalent to $\Omega^\infty{\bf
Der}_R(X,\Sigma I)$, see Remark \ref{inf}. Since according to our
assumption the space \[\xymatrix{\hofib
\Map_{\mathcal{C}^{ass}_R}(X,B)\ar[r]&
\Map_{\mathcal{C}^{ass}_R}(X,A)}\] is nonempty the images of the
lower and right arrows in (\ref{usd}) coincide and therefore
\begin{align*}\hofib \Map_{\mathcal{C}^{ass}_R}(X,B)\ra
\Map_{\mathcal{C}^{ass}_R}(X,A)&\simeq \Omega (\Omega^\infty{\bf
Der}_R(X, \Sigma I))\\&\simeq \Omega^\infty{\bf
Der}_R(X,I).\end{align*}  This finishes the proof of Theorem
\ref{obs}.

A large supply of topological singular extensions can be obtain by
taking Postnikov stages of connective $R$-algebras of commutative
$R$-algebras. We have the following
\begin{theorem}\label{Postnikov} Let $A$ be a connective $R$-algebra where $R$ is a
connective commutative $S$-algebra. Then there exists a tower of
$R$-algebras under $A$:
\[\xymatrix{&&A\ar_{f_0}[dll]\ar^{f_1}[dl]\ar_{f_n}[dr]\ar^{f_{n+1}}[drr]\\ A_0&
A_1\ar[l]&\ldots\ar[l]&A_n\ar[l]& A_{n+1}\ar[l]&\ldots\ar[l]}\]
such that \begin{itemize} \item $\pi_kA_n=0$ for $k>n$; \item the
map $f_n:A\lra A_n$ induces an isomorphism $\pi_kA\lra \pi_kA_n$
for $0\leq k\leq n$; \item the homotopy fibre sequences
$H\pi_{n+1}A\lra A_{n+1}\lra A_{n}$ are topological singular
extensions associated with topological derivations $k_n\in
Der^1_R(A_n, H\pi_{n+1}A)$.
\end{itemize}
\end{theorem}
In the case when $A$ is a connective commutative $R$-algebra we
have the precise analogue of the theorem above: there is a tower
$\{A_n\}$ of \emph{commutative} $R$-algebras under $A$ and the
homotopy fibre sequences $H\pi_{n+1}A\lra A_{n+1}\lra A_{n}$ are
topological singular extensions associated with \emph{commutative}
derivations $k_n\in TAQ^1_R(A_n, H\pi_{n+1}A)$.

Theorem \ref{Postnikov} is proved by induction on $n$. By glueing
cells in the category of $R$-algebras or commutative $R$-algebras
we construct a map $A\lra A_0=H\pi_0A$. Suppose that the $n$th
stage $A_n$ for a commutative $R$-algebra $A_n$ is already
constructed. It is not hard to see using the spectral sequence
(\ref{sp1}) that the lowest nonzero homotopy group of
$\Omega^{comm}_{A_n/A}$ is in dimension $n+1$ and is equal to
$\pi_{n+1}A$. We then construct a map
$k_n:\Omega^{comm}_{A_n/A}\lra H\pi_{n+1}A$ by attaching cells to
$\Omega^{comm}_{A_n/A}$ to kill its higher homotopy groups. The
map $k_n$ is the required commutative derivation. In the
associative context we use the $A$-bimodule of relative
differentials $\Omega^{ass}_{A_n/A}$ and proceed similarly.

Theorem \ref{Postnikov} provides a large supply of $R$-algebras or
commutative $R$-algebras. For example Postnikov stages of $MU$ or
$kU$ are commutative $S$-algebras. Other examples of commutative
$S$-algebras include the Morava $E$-theory spectra $E_n$ and
Eilenberg-MacLane spectra $Hk$ for a commutative ring $k$. It is
generally very hard to prove that a spectrum possesses a
commutative $S$-algebra structure unless there is a geometric
reason for it. We will see later on that the situation with
associative structures is better; a large class of
complex-oriented spectra could be given structures of
$MU$-algebras and, therefore, $S$-algebras.
\section{Stabilization}
In this section we give an interpretation of topological
derivation and topological Andr\'e-Quillen cohomology as
`stabilizations' of a certain forgetful functors. This
interpretation in the commutative case is due to Basterra-McCarthy
\cite{BM} and Basterra-Mandell (private communication).

Let $F:\mathcal{C}\lra \mathcal{D}$ be a continuous functor
between pointed topological model categories. We assume that $F$
is a homotopy functor which means that it preserves homotopy
equivalences. We also assume that $F$ is reduced, that is it takes
the initial object of $\mathcal{C}$ into a contractible object in
$\mathcal{D}$ (i.e. homotopy equivalent to the initial object in
$\mathcal{D}$). For an object $X$ in $\mathcal{D}$ we will denote
by $\Omega X$ the cotensor of $X$ and the pointed circle $S^1$. In
addition the tensor of $?$ with a {pointed} topological space $X$
will be denoted by $X\hat{\otimes}?$.

We have the following pushout diagram of pointed spaces:
\begin{equation}\label{push}\xymatrix{S^0\ar[r]\ar[d]&I\ar[d]\\I\ar[r]&S^{1}}
\end{equation} where $I$ is the unit interval. Tensoring it with $X$ and
applying the functor $F$ we get the diagram
\[\xymatrix{F(X)=F(S^0\hat{\otimes}X)\ar[r]\ar[d]&F(I\hat{\otimes}X)
\ar[d]\\F(I\hat{\otimes}X)\ar[r]&F(S^{1}\hat{\otimes}X)}\] Define
the functor $X\lra TF(X)$ by the requirement that the diagram
\[\xymatrix{TF(X)\ar[r]\ar[d]&F(I\hat{\otimes}X)
\ar[d]\\F(I\hat{\otimes}X)\ar[r]&F(S^{1}\hat{\otimes}X)}\] be a
homotopy pullback in $\mathcal{D}$. There results a natural
transformation $\xi:F(X)\lra TF(X)$.
\begin{defi} The stabilization of $F$ is the functor $F^{st}$
defined as
\[\xymatrix{F^{st}(X)=\hocolim(F(X)\ar[r]&TF(X)\ar[r]&TTF(X)\ar[r]&\ldots)}\]
where the colimit is taken over the maps $T^i(\xi): T^i(F(X))\lra
T^i(TF(X))$.
\end{defi}

Note that since $I\hat{\otimes}X$ is contractible the object
$TF(X)$ is homotopy equvalent to $\Omega F(S^1\hat{\otimes} X)$.
Furthermore, it is easy to see that there is a homotopy
equivalence
\[T^nFX\simeq \Omega^nF(S^n\hat{\otimes}X)\] where $\Omega^n$ is
the $n$th iterate of $\Omega$.
\begin{rem} The stabilization (or linearization) of a homotopy
functor was introduced by Goodwillie in \cite{Goo1} and
constitutes the first layer of its \emph{Taylor tower},
\cite{Goo3}. It would be interesting to investigate the higher
layers of this tower in our context.\end{rem} We are now ready to
give the interpretation of  Der and TAQ as stabilizations of
certain forgetful functors. We start with the associative case.
Let $A$ be a cofibrant $R$-algebra and consider the category
$\mathcal{C}^{ass}_{A/A}$ of $A$-algebras over $A$. Recall that
$\mathcal{C}^{ass}_{A/A}$ consists of $R$-algebras $B$ supplied
with two maps $\eta_B:B\lra A$ and $\zeta_B:A\lra B$ such that
$\eta_A\circ\zeta_A=id_A$. The category $\mathcal{C}^{ass}_{A/A}$
is a topological model category. Let us first describe its
cotensors. Let $\tilde{B}$ be a fibrant approximation of $B$ as an
object in $\mathcal{C}^{ass}_{A/A}$. Then
$\eta_{\tilde{B}}:\tilde{B}\lra A$ is a fibration of $R$-algebras.
Denote by $I(B)$  the fibre of $\eta_{\tilde{B}}$ so that there is
a homotopy fibre sequence of $B$-bimodules:
\[\xymatrix{I(B)\ar[r]&B\ar^{\eta_B}[r]&A}.\] Consider the $R$-module
$A\vee I(B)$. It is clearly an object in $\mathcal{C}^{ass}_{A/A}$
and, moreover, the obvious map $A\vee I(B)\lra B$ is a weak
equivalence. Then for a based space $X$ the $R$-algebra $B^X$, the
cotensor of $B$ and $X$ in $\mathcal{C}^{ass}_{A/A}$ is weakly
equivalent to $A\vee I(B)^X$ where $I(B)^X=F_S(\Sigma^{\infty}X,
I(B))$ is the usual cotensor (function spectrum) of $X$ and the
$R$-module $I(B)$.

Now we will describe tensors in $\mathcal{C}^{ass}_{A/A}$. Let $B$
be a cofibrant object in $\mathcal{C}^{ass}_{A/A}$ and $X$ be a
based space that consists of $n+1$ points, one of them being the
base point. Then the tensor $X\hat{\otimes}B$ is just the n-fold
free product of $B$ over $A$: $B^{\coprod_An}$. For a general
based space $X$ we replace it with the simplicial set $X_*$ whose
realization is equivalent to $X$, and tensor it component-wise
with $B$. The realization of the obtained simplicial object in
$\mathcal{C}^{ass}_{A/A}$ will be homotopy equivalent to
$X\hat{\otimes}B$. For example if we take $X=S^1$, the circle and
by $X_*$ its simplicial model having two nondegenerate simplices
in dimensions $0$ and $1$ then the resulting simplicial object
$X_*\hat{\otimes}B$ is just the bar-construction $\beta_*(A,B,A)$.
We have
\[\beta_i(A,B,A)=B^{ \coprod_Ai}.\] The face maps are induced in the usual way
by the canonical folding map $B\coprod_AB\lra B$ and the
degeneracies are induced by the map $\zeta_B$.

Recall that  the functor $\Omega^{ass}_{B/A}$ is defined from the
fibre sequence
\[\xymatrix{\Omega^{ass}_{B/A}\ar[r]&B\wedge_{A} B\ar^-{m^\prime}[r]& B}.\] We will
assume that $\Omega^{ass}_{B/A}$ is a cofibrant $B$-bimodule or
else replace it by its cofibrant approximation. By definition,
\[\Omega^{ass}_{A/A}(B):=A\wedge_B\Omega^{ass}_{B/A}\wedge_BA.\] Both $I$ and
$\Omega^{ass}_{A/A}(B)$ considered as functors of $B$ take their
values in the category of $A$-bimodules.
\begin{theorem}\label{stab} Let $A$ and $B$ be connective $R$-algebras. Then
$I^{st}(B)\simeq\Omega^{ass}_{A/A}(B)$\end{theorem} To see this
consider the `universal derivation' map $\xymatrix{B\ar[r]&A\vee
\Omega^{ass}_{A/A}(B)}$ and the composition
\[f_B:\xymatrix{I(B)\ar[r]&B\ar[r]&A\vee\Omega^{ass}_{A/A}(B)\ar[r]&
\Omega^{ass}_{A/A}(B)}.\]
It is easy to see that there is a weak equivalence of  $R$-modules
\[\Omega^{ass}_{A/A}(B)\simeq A\wedge_BI(B)\] and the map $f_B$ could be
represented as \[\xymatrix{I(B)\cong
B\wedge_BI(B)\ar^-{\eta_B\wedge id}[r]&A\wedge_BI(B).}\]Note that
if the $\eta_B:B\lra A$ is $n$-connected then the map $f_B$ is
$2n+1$-connected. Further if $I(B)$ is $m$-connected then
$I(S^n\hat{\otimes}B)$ is $n+m$-connected and the map
\[\xymatrix{f_{S^n\hat{\otimes}B}:I(S^n\hat{\otimes}B)\ar[r]&
\Omega_{A/A}^{ass}(S^n\hat{\otimes}B)}\] is $2(m+n)+1$-connected.
It follows that the map
\[\xymatrix{\Omega^nf_{S^n\hat{\otimes}B}:\Omega^nI(S^n\hat{\otimes}B)\ar[r]&
\Omega^n\Omega_{A/A}^{ass}(S^n\hat{\otimes}B)}\] is
$2m+1$-connected and therefore
$I^{st}(B)=[\Omega^{ass}_{A/A}]^{st}(B)$. The results of Section 2
imply that the functor $B\rightarrow [\Omega^{ass}_{A/A}]^{st}(B)$
preserves tensors and therefore
\[\Omega^n\Omega^{ass}_{A/A}(S^n\hat{\otimes}B)\simeq \Omega^n\Sigma^n
\Omega^{ass}_{A/A}(B)\simeq \Omega^{ass}_{A/A}(B).\] which implies
that ${[\Omega^{ass}_{A/A}}]^{st}(B)\simeq \Omega^{ass}_{A/A}(B)$.
This finishes our sketch of the proof of Theorem \ref{stab}.

What is the relation of the functor $\Omega^{ass}_{A/A}(B)$ to the
topological derivations of $A$? Let $B=A\coprod A$. Define
$\zeta_B$ to be the map $\xymatrix{A=A\coprod R\ar^-{id\coprod
1}[r]&A\coprod A}$ and $\eta_B:A\coprod A\lra A$ to be the folding
map. Then $B$ becomes an object in $\mathcal{C}^{ass}_{A/A}$.
Moreover we have a natural isomorphism for any $C\in
\mathcal{C}^{ass}_{A/A}$:
\[h\mathcal{C}^{ass}_{A/A}(B,C)\cong h\mathcal{C}_{R/A}(A,C).\]

Therefore for an $A$-bimodule $M$ we have
\begin{align*}Der^0_R(A,M)=h\mathcal{C}_{R/A}(A,A\vee M) \cong&
h\mathcal{C}^{ass}_{A/A}(B,A\vee M) \\ \cong&
h\mathcal{M}_{A-A}(\Omega^{ass}_{A/A},M)\\ \cong&
h\mathcal{M}_{A-A}(I^{st}(B),M).\end{align*} We can also obtain in
the usual way the enriched version of the above equivalence:
\[{\bf Der}_R(A,A\vee M)\simeq
F_{A\wedge A^{op}}(I^{st}(A\coprod A),M).\] Theorem \ref{stab} has
an analogue in the context of commutative $R$-algebras. The
corresponding result was first established by M.Basterra and
M.Mandell, cf, also \cite{BM} for TAQ  and our proof was modeled
on theirs. One should also mention that the original definition of
TAQ due to I. Kriz \cite{Kriz} was in terms of stabilization. In
contrast with the associative case this interpretation is
extremely helpful in concrete calculations.

Let us now describe the result of Basterra and Mandell. $A$ is now
a commutative cofibrant $R$-algebra. Consider the category
$\mathcal{C}^{comm}_{A/A}$ of commutative $A$-algebras over $A$.
Then $I(B)\in \mathcal{M}_A$ is defined from the  fibre sequence
\[\xymatrix{I(B)\ar[r]&\tilde{B}\ar^{\eta_{\tilde{B}}}[r]&A}\] where
$\tilde{B}$ is the fibrant replacement of $B$ in
$\mathcal{C}^{comm}_{A/A}$. Recall the definition of the
\emph{indecomposables functor} $B\rightarrow Q(B)$ from Section 1.
It turns out that if $A$ and $B$ are connected then $I^{st}(B)$,
the stabilization of the functor $I$ evaluated at $B$, is weakly
equivalent to $Q\circ I(B)$. The proof is very similar to the one
given in the associative case. To show that the map

\[\xymatrix{I(S^n\hat{\otimes}B)\ar[r]&Q\circ I(S^n\hat{\otimes}B)}\]
is sufficiently highly connected one uses Proposition \ref{har}.
To relate stabilization with TAQ consider $A\wedge A$ as an
$A$-algebra over $A$. The augmentation $A\wedge A\lra A$ is just
the multiplication map and the $A$-algebra structure on $A\wedge
A$ is defined via the action of $A$ on the left smash factor. Then
for an $A$-module $M$ we have the following isomorphisms in
$h\mathcal{M}_A$:
\[{\bf TAQ}_R(A,M)\simeq F_A(Q\circ I(A\wedge A),M)\simeq
F_A(I^{st}(A\wedge A),M).\] In fact the author was informed by
M.Basterra that the last isomorphism holds even without the
hypothesis that $A$ is connective. That suggests that Theorem
\ref{stab} might also hold without any connectivity assumptions.

\section{Calculations} Up until now we considered
only formal properties of TAQ and topological derivations and saw
that they are quite similar. The difference appears in their
calculational aspect. Topological Andr\'e-Qullen cohomology is far
harder to compute and presently very few explicit calculations
have been carried out. By contrast, topological derivations are
relatively approachable. In most cases their computation reduces
to that of topological Hochschild cohomology thanks to the
homotopy fibre sequence (\ref{hoch}). There is an extensive
literature dedicated to the computation of THH of different
Eilenberg-MacLane spectra of which we mention \cite{Bok},
\cite{Fra}, \cite{Mad}. Topological Hochschild (co)homology of
some spectra other than Eilenberg-MacLane are computed in
\cite{MS} and \cite{BaL}.  We do not intend to give a detailed
review of this literature and restrict ourselves with giving a few
examples related to obstruction theory. In this section the symbol
$`\wedge$' stands for $`\wedge_S$'.

Our first example is the computation of topological derivations of
$H\mathbb{F}_p$.
\begin{theorem}\label{bork} \begin{enumerate}\item There is an isomorphism of
graded $\mathbb{F}_p$-vector spaces
\[THH_S^*(\mathbb{F}_p,\mathbb{F}_p)\cong
(\mathbb{F}_p[x])^*=Hom_{\mathbb{F}_p}(\mathbb{F}_p[x],\mathbb{F}_p).\]
where $x$ has degree $2$. \item There is an isomorphism of graded
$\mathbb{F}_p$-vector spaces
\[Der_S^{*-1}(\mathbb{F}_p,\mathbb{F}_p)\cong
(\mathbb{F}_p[x])^*/(\mathbb{F}_p)\]\end{enumerate}\end{theorem}
Note that the statement (2) of the theorem follows from statement
(1) because the homotopy fibre sequence (\ref{hoch}) in our case
splits giving the short exact sequence
\[\xymatrix{\mathbb{F}_p\ar[r]&THH_S^*(\mathbb{F}_p,\mathbb{F}_p)\ar[r]&
Der_S^{*-1}(\mathbb{F}_p,\mathbb{F}_p)}.\] To obtain the claim
about $THH_S^*(\mathbb{F}_p,\mathbb{F}_p)$ it is enough to compute
topological Hochschild \emph{homology}
$THH^S_*(\mathbb{F}_p,\mathbb{F}_p)$ because
\[THH_S^*(\mathbb{F}_p,\mathbb{F}_p)\cong
Hom_{\mathbb{F}_p}(THH^S_*(\mathbb{F}_p,\mathbb{F}_p),\mathbb{F}_p),\]
see Remark \ref{dual}. The computation of
$THH^S_*(\mathbb{F}_p,\mathbb{F}_p)$ is due to B\"okstedt,
\cite{Bok}. Briefly, one uses the spectral sequence
\[E^2_{**}=Tor_{{H\mathbb{F}_p\wedge
H\mathbb{F}_p}_*}(\mathbb{F}_p,\mathbb{F}_p)=
Tor_{\mathcal{A}_*^p}(\mathbb{F}_p,\mathbb{F}_p)\Rightarrow
THH^S_*(\mathbb{F}_p,\mathbb{F}_p)\] where
$\mathcal{A}_*^p={H\mathbb{F}_p\wedge H\mathbb{F}_p}_*$ is the
dual Steenrod algebra mod $p$. The $E_2$-term of this spectral
sequence is easy to compute. It turns out that for $p=2$ there are
no further differntials and for an odd $p$ one uses the
Dyer-Lashof operations to compute the only nontrivial differential
$d^{p-1}$. The same operations also solve the extension problem.
The final result is that, regardless of $p$,
\begin{equation}\label{bo}THH^S_*(H\mathbb{F}_p,H\mathbb{F}_p)\cong
\mathbb{F}_p[x]\end{equation} with $x$ in homological degree $2$.
This finishes our sketch of Theorem \ref{bork}.
\begin{rem}\label{help} In fact, B\"okstedt's calculation shows that the
isomorphism of (\ref{bo}) is multiplicative. This fact will be
used in the computation of
$TAQ^*_S(H\mathbb{F}_p,H\mathbb{F}_p)$.\end{rem}

Our next example is the calculation of ${ Der
}^*_{MU}(H\mathbb{F}_p,H\mathbb{F}_p)$ where $MU$ is the complex
cobordism spectrum. Recall that the coefficient ring of $MU$ is
the polynomial algebra:\[MU_*=\mathbb{Z}[x_1,x_2,\ldots]\] where
$|x_i|=2i$. Since $MU$ is known to be a commutative connective
$S$-algebra the spectrum $H\mathbb{F}_p$ is naturally an
$MU$-algebra and it makes sense to consider topological
derivations and topological Hochschild cohomology of
$H\mathbb{F}_p$ as an $MU$-algebra. The relevant result is
\begin{theorem}\label{MU}\begin{enumerate}\item There is an isomorphism of
graded $\mathbb{F}_p$-vector spaces
\[THH^*_{MU}(H\mathbb{F}_p,H\mathbb{F}_p)=\mathbb{F}_p[y_1,y_2,\ldots]\]
where $y_i$ has cohomological degree $2i$. \item There is an
isomorphism of graded $\mathbb{F}_p$-vector
spaces\[Der^{*-1}_{MU}(H\mathbb{F}_p,H\mathbb{F}_p)=
\mathbb{F}_p[y_1,y_2,\ldots]/(\mathbb{F}_p).\]\end{enumerate}\end{theorem}
Again the second claim is a consequence of the first one. To
verify the first claim we need to compute
$\pi_*H\mathbb{F}_p\wedge_{MU}H\mathbb{F}_p$. Using the spectral
sequence
\begin{align*}Tor_{**}^{MU_*}(\mathbb{F}_p,\mathbb{F}_p)&=
Tor_{**}^{\mathbb{Z}[x_1,x_2,\ldots]}(\mathbb{F}_p,\mathbb{F}_p)\\
&\Rightarrow
\pi_*H\mathbb{F}_p\wedge_{MU}H\mathbb{F}_p\end{align*} we see that
the ring $\pi_*H\mathbb{F}_p\wedge_{MU}H\mathbb{F}_p$ is
isomorphic to $\Lambda_{\mathbb{F}_p}(z_1,z_2,\ldots)$, the
exterior algebra over $\mathbb{F}_p$ on generators $z_i$ where
$|z_i|=2i-1$.

Finally from the spectral sequence
\begin{align*}
Ext^{**}_{\pi_*H\mathbb{F}_p\wedge_{MU}H\mathbb{F}_p}(\mathbb{F}_p,\mathbb{F}_p)&=
Ext^{**}_{\Lambda_{\mathbb{F}_p}(z_1,z_2,\ldots)}(\mathbb{F}_p,\mathbb{F}_p)\\
&=\mathbb{F}_p[y_1,y_2,\ldots]\Rightarrow
THH^*_{MU}(H\mathbb{F}_p,H\mathbb{F}_p)\end{align*} we obtain the
desired isomorphism.

Next we would like to discuss the forgetful maps
\[\xymatrix{l_{MU}^*:Der^*_{MU}(H\mathbb{F}_p,H\mathbb{F}_p)\ar[r]&
[H\mathbb{F}_p,H\mathbb{F}_p]^*_{MU-mod}}\] and
\[\xymatrix{l_S^*:Der^*_{S}(H\mathbb{F}_p,H\mathbb{F}_p)\ar[r]&
[H\mathbb{F}_p,H\mathbb{F}_p]^*_{S-mod}=\mathcal{A}^*_p}.\] We
already saw that the algebra of cooperations in $MU$-theory
$\pi_*H\mathbb{F}_p\wedge_{MU}H\mathbb{F}_p$ is an exterior
algebra on generators $z_i$. Just as the usual dual Steenrod
algebra it is a Hopf algebra and it is easy to see that $z_i$'s
are primitive elements in it. Therefore  the $MU$-Steenrod algebra
$[H\mathbb{F}_p,H\mathbb{F}_p]^*_{MU-mod}$ is itself an exterior
algebra $\Lambda_{\mathbb{F}_p}(z^*_1,z^*_2,\ldots)$ where
$z^*_i$'s are dual to $z_i$.

Then the map $l_{MU}^*$ sends the derivation $y_i\in
Der^{2i-1}_{MU}(H\mathbb{F}_p,H\mathbb{F}_p)$ to the operation
$z^*_i\in [H\mathbb{F}_p,H\mathbb{F}_p]^{2i-1}_{MU-mod}$.

Any map of $MU$-modules $H\mathbb{F}_p\lra\Sigma^*H\mathbb{F}_p$
can also be considered as a map of $S$-modules. In other words
there is a map from the $MU$-Steenrod algebra into
$\mathcal{A}^*_p$:
\[\xymatrix{[H\mathbb{F}_p,H\mathbb{F}_p]^*_{MU-mod}\ar[r]&
[H\mathbb{F}_p,H\mathbb{F}_p]^*_{S-mod}=\mathcal{A}^*_p}.\] The
latter map is a map of Hopf algebras and should therefore respect
primitive elements. We know that the odd degree primitive elements
are the Milnor Boksteins $Q_i\in\mathcal{A}^*_p$. Since
$|Q_i|=2p^i-1$ only the primitives $z^*_{2p^i-1}$ could have a
nonzero image in $\mathcal{A}^*_p$. To see that $z^*_{2p^i-1}$
does indeed have a nonzero image notice that it is the first
nonzero Postnikov invariant of the connective Morava $K$-theory
spectrum $k(n)$. This invariant is nontrivial both in the category
of $S$-modules and $MU$-modules which shows that the image of
$z^*_{2p^i-1}$ is $Q_i$ up to an invertible scalar factor.

This allows one to calculate the image of the map $l_S^*$. On the
one hand its image should be contained in the subspace spanned by
the primitive odd degree elements in $\mathcal{A}^*_p$, that is,
the elements $Q_i$. On the other hand any operation $Q_i$ can be
lifted to a topological derivation in
$Der^*_S(H\mathbb{F}_p,H\mathbb{F}_p)$ since it even lifts to an
$MU$-derivation in $Der^*_{MU}(H\mathbb{F}_p,H\mathbb{F}_p)$.
Therefore the image of $l_S$ is the whole subspace of
$\mathcal{A}^*_p$ spanned by $Q_i$'s.

We now turn to the calculation, due to M.Basterra and M.Mandell,
of $TAQ^*_S(H\mathbb{F}_p,H\mathbb{F}_p)$. We have:
\begin{align*}TAQ^*_S(H\mathbb{F}_p,H\mathbb{F}_p)&\simeq [\hocolim_{n\rightarrow
\infty}\Omega^nI(S^n\hat{\otimes}(H\mathbb{Z}\wedge
H\mathbb{F}_p)),H\mathbb{F}_p ]^*_{H\mathbb{F}_p-mod}\\ &\simeq
[\hocolim_{n\rightarrow
\infty}\Sigma^{-n}I(S^{n-2}\hat{\otimes}[S^2\hat{\otimes}(H\mathbb{F}_p\wedge
H\mathbb{F}_p)]),H\mathbb{F}_p]^*_{H\mathbb{F}_p-mod}
\\ &\simeq
[\Sigma^{-2}\hocolim_{n\rightarrow \infty}
\Sigma^{-n+2}I(S^{n-2}\hat{\otimes}[S^2\hat{\otimes}(H\mathbb{F}_p\wedge
H\mathbb{F}_p)]),H\mathbb{F}_p]^*_{H\mathbb{F}_p-mod}
\\&\simeq
\Sigma^{-2}TAQ^*_{H\mathbb{F}_p}(S^2\hat{\otimes}(H\mathbb{F}_p\wedge
H\mathbb{F}_p),H\mathbb{F}_p).
\end{align*}
\begin{lem}\label{MB} There is a weak equivalence of commutative
$H\mathbb{F}_p$-algebras
\begin{equation}\label{equiv}S^2\hat{\otimes}(H\mathbb{F}_p\wedge H\mathbb{F}_p)\cong
H\mathbb{F}_p\vee\Sigma^3H\mathbb{F}_p.\end{equation} Here the
right hand side of (\ref{equiv}) is the square-zero extension of
$H\mathbb{F}_p$ by $\Sigma^3H\mathbb{F}_p$. \end{lem} To see this
note first that \[S^1\hat{\otimes}(H\mathbb{F}_p\wedge
H\mathbb{F}_p)\simeq {\bf THH}^S(H\mathbb{F}_p,H\mathbb{F}_p).\]
Therefore \[\pi_*S^1\hat{\otimes}(H\mathbb{F}_p\wedge
H\mathbb{F}_p)=\mathbb{F}_p[x]\] with $|x|=2$, see Remark
\ref{help}.  Denote the commutative $S$-algebra
$S^1\hat{\otimes}(H\mathbb{F}_p\wedge H\mathbb{F}_p)$ by $X$. We
have the following weak equivalences of commutative $S$-algebras:
\[ S^2\hat{\otimes}(H\mathbb{F}_p\wedge H\mathbb{F}_p)\simeq
S^1\hat{\otimes}X\simeq H\mathbb{F}_p\wedge _XH\mathbb{F}_p.\]
From the spectral sequence
\[Tor_{**}^{X_*}(\mathbb{F}_p,\mathbb{F}_p)=
\Lambda_{\mathbb{F}_p}(y)\Rightarrow \pi_*H\mathbb{F}_p\wedge
_XH\mathbb{F}_p\] we see that
\[\pi_*H\mathbb{F}_p\wedge
_XH\mathbb{F}_p=\pi_*S^2\hat{\otimes}(H\mathbb{F}_p\wedge
H\mathbb{F}_p)=\Lambda_{\mathbb{F}_p}(y).\] Here $y$ has degree
$3$. It is not hard to see, using Proposition \ref{har} that any
augmented $H\mathbb{F}_p$-algebra whose coefficient ring is an
exterior algebra on one generator $y$ in positive degree is in
fact weakly equivalent to the square-zero extension
$H\mathbb{F}_p\vee \Sigma^{|y|}H\mathbb{F}_p.$ This finishes our
sketch of the proof of Lemma \ref{MB}.

It follows that
\[TAQ^*_S(H\mathbb{F}_p,H\mathbb{F}_p)=
\Sigma^{-2}TAQ^*_{H\mathbb{F}_p}(H\mathbb{F}_p\vee \Sigma^3
H\mathbb{F}_p,H\mathbb{F}_p).\] So we need to compute the
topological Andr\`e-Quillen cohomology of the `Eilenberg-MacLane
object' $H\mathbb{F}_p\vee \Sigma^3H\mathbb{F}_p$ in the category
of commutative augmented $H\mathbb{F}_p$-algebras. There is an
obvious element (`fundamental class') $x\in
TAQ^{3}_{H\mathbb{F}_p}(H\mathbb{F}_p\vee \Sigma^3
H\mathbb{F}_p,H\mathbb{F}_p)$ corresponding to the identity map on
$H\mathbb{F}_p\vee \Sigma^3 H\mathbb{F}_p$. Further one introduces
the action of Steenrod operations on TAQ. Let $p=2$. Then there
exist operations
\[\xymatrix{Sq^i:TAQ^q_{H\mathbb{F}_p}(-, H\mathbb{F}_p)\ar[r]&
TAQ^{q+i}_{H\mathbb{F}_p}(-,H\mathbb{F}_p)}\] defined for $q<i$
and satisfying the usual Adem relations. Applying $Sq^i$ to
$y=\Sigma^{-2}x\in
TAQ^{1}_{H\mathbb{F}_p}(H\mathbb{F}_p,H\mathbb{F}_p)$ one obtains
new elements in
$TAQ^*_{H\mathbb{F}_p}(H\mathbb{F}_p,H\mathbb{F}_p)$. The final
result is (for $p=2$): \begin{theorem}
$TAQ^*_S(H\mathbb{F}_p,H\mathbb{F}_p)$ is spanned by the elements
$Sq^{s_1}Sq^{s_2}\ldots Sq^{s_r}y$ where the sequence
$(s_1,s_2,\ldots,s_r)$ is Steenrod admissible and $s_r>3$.
\end{theorem}
There is a similar result for odd primes which uses the odd
primary operations in TAQ instead of Steenrod squares.

To conclude this section note that the forgetful map
\[\xymatrix{TAQ^*_S(H\mathbb{F}_p,H\mathbb{F}_p)\ar[r]&
[H\mathbb{F}_p,H\mathbb{F}_p]^*_{S-mod}=\mathcal{A}_p^*}\] is not
as interesting as in the associative case. It is easy to see that
its image contains  (multiples of) the Bokstein homomorphism
$\beta$ corresponding to the topological singular extension
\[\xymatrix{H\mathbb{F}_p\ar[r]& H\mathbb{F}_{p^2}\lra
H\mathbb{F}_p}\] and nothing else. Indeed, if it contained, e.g.
the operation $Q_i$ for some $i>0$ then the corresponding singular
extension $H\mathbb{F}_p\lra ?\lra \Sigma^{2(p^i-1)}H\mathbb{F}_p$
would realize the $2$-stage Postnikov tower of $k(n)$ as a
commutative $S$-algebra which is impossible.
\section{Existence of $S$-algebra structures}
Topological Andr\`e-Quillen cohomology was originally introduced
in \cite{Kriz} with the purpose of proving that $BP$, the
Brown-Peterson spectrum supports a structure of an $E_\infty$-ring
spectra (or, equivalently, a commutative $S$-algebra structure), a
long-standing problem posed by P. May. Unfortunately, this problem
is still open, however the problem concerning \emph{associative}
structures turned out to be much more manageable. The first result
in this direction is the seminal paper of A. Robinson \cite{Rob1}
where Morava $K$-theory spectra at odd primes were proved to
support an $A_\infty$-structure. There are some difficulties in
interpreting Robinson's proof within our context, particularly the
existence of a unit seems problematic. However working in the
category of $MU$-modules it is possible to obtain considerably
stronger results. Namely, in \cite{L2} and \cite{Goe} it was
proved that a broad class of complex-oriented cohomology theories
admit structures of $S$-algebras. These structures are typically
non-unique and in the last section we say something about the
corresponding moduli spaces. Here we describe the relevant results
following \cite{L2}.

Let us assume that the coefficient ring $R_*$ of our `ground'
commutative $S$-algebra $R$ is concentrated in even degrees. For
an element $x\in R_*$ we will denote by $R/x$ the
 cofibre of the map $R\stackrel{x}{\rightarrow}R$.
Let $I_\ast$ be a graded ideal generated by a (possibly infinite)
regular sequence of homogeneous elements $(u_1,u_2,\ldots)\in
R_\ast$. We assume in addition that each $u_k$ is a nonzero
divisor in the ring $R_*$. Then we can form the $R$-module $R/I$
as the infinite smash product of $R/u_k$. It is known by work of
Strickland \cite{Str} that there is a structure of an $R$-ring
spectrum on $R/I$. Clearly the coefficient ring of $R/I$ is
isomorphic to $R_\ast/I_\ast$ where $R_\ast/I_\ast$ is understood
to be the direct limit of $R_\ast/(u_k,u_2,\ldots,u_k)$. Let us
denote the $R$-algebra $R/I$ by $A$. Our standing assumption is
that $A$ has a structure of an $R$-algebra (i.e. strictly
associative). This may seem a rather strong condition but, as we
see shortly such a situation is quite typical. In fact P.Goerss
proved in \cite{Goe} that any spectrum obtained by killing a
regular sequence in $MU$, the complex cobordism $S$-algebra, has a
structure of an $MU$-algebra.

The construction we are about to describe allows one to construct
new $R$-algebras by `adjoining' the indeterminates $u_k$ to the
$R$-algebra $A$.
 Let us introduce
the notation $A(u^l_k)_*$ for the $R_\ast$-algebra
$$\lim_{n\rightarrow \infty}R_\ast/(u_1,u_2,
\ldots,u_{k-1},u_k^l,u_{k+1},\ldots, u_n).$$
 For each $l$ reduction modulo
$u_k^{l}$ determines a map of $R_\ast$-algebras
\[{A(u^{l+1}_k)_*}\rightarrow {A(u^l_k)_*}. \] Now we can formulate
our main theorem in this section.
\begin{theorem} \label{qas}
For each $l$ there exist $R$-algebras $A(u^l_k)$ with coefficient
rings ${A(u^l_k)_*}$ and $R$-algebra maps
\[R_{l,k}:A(u^{l+1}_k)\rightarrow A(u^l_k)\] which give the
reductions mod $u_k^{l}$ on the level of coefficient rings.
\end{theorem}
Let us explain the idea of the proof. One uses induction on $l$.
The case $l=1$ is just our original assumption that $A$ is an
$R$-algebra. Suppose that the $R$-algebras $A(u^l_k)$ with the
required properties were constructed for $l\leq i$. One then shows
that $R/u^{i+1}_k$ possesses an $R$-ring structure and there
exists a homotopy cofibre sequence of $R$-modules
\begin{equation}\label{bk1}
\xymatrix{R/u^{i+1}_k\ar[r]&R/u^{i}_k\ar[r]&\Sigma^{di+1}R/u_k}\end{equation}
where the first arrow is an $R$-ring spectrum map realizing the
reduction mod $u^{i}$ in homotopy and the second arrow is an
appropriate Bokstein operation. Taking the smash product of
(\ref{bk1}) with $R/u_1\wedge R/u_2\wedge\ldots
 \wedge R/u_{k-1}\wedge R/u_{k+1}\wedge\ldots$ we get the following
 homotopy cofibre sequence
of $R$-modules:
\begin{equation}\label{bk}\xymatrix{A(u^{i+1}_k)\ar[r]& A(u^i_k)
\ar[r]&\Sigma^{di+1}A.}\end{equation} Finally one proves that the
sequence (\ref{bk}) could be improved to a topological singular
extension by computing $Der^*_R(A(u^i_k), A)$ and analyzing the
forgetful map \[\xymatrix{Der^*_R(A(u^i_k), A)\ar[r]&[A(u^i_k),
A]^*_{R-mod}}.\]

Theorem \ref{qas} leads to a host of examples of $S$-algebras. Let
$R=MU_{(p)}$. the complex cobordism spectrum localized at $p$ and
$I_*=(p, x_1,x_2,\ldots)$ be the maximal ideal in
$MU_{(p)*}=\mathbb{Z}_{(p)}[x_1,x_2,\ldots]$. Then $MU/I\simeq
H\mathbb{F}_p$ is an $MU$-algebra. For any polynomial generator
$x_k\in MU_{(p)}$ and any $i$ we can construct an
$MU_{(p)}$-algebra realizing in homotopy the $MU_{(p)*}$-module
$\mathbb{F}_p[x_k^i]$. Taking the homotopy inverse limit of these
we get an $MU_{(p)}$-algebra whose coefficient ring is
$\mathbb{F}_p[x_k]$. Iterating this procedure we can realize
topologically the $MU_{(p)*}$-algebra
$\mathbb{F}_p[x_{k_1},x_{k_2},\ldots]$ for any sequence (possibly
infinite) of polynomial generators $x_{k_1},x_{k_2},\ldots.$ We
could also set $I_*=(x_1,x_2,\ldots)$ to obtain integral versions
of these $MU_{(p)}$-algebras.

Recall that the Brown-Peterson spectrum $BP$ is obtained from
$MU_{(p)}$ by killing all polynomial generators in $MU_{(p)*}$
except for $v_i=x_{2(p^i-1)}$ provided that we use Hazewinkel
generators for $MU_{(p)*}$. Therefore
$BP_\ast=\mathbb{Z}_{(p)}[v_1,v_2,\ldots]$. It follows that $BP$
possesses an $MU_{(p)}$-algebra structure. We define \[BP\langle
n\rangle=BP/(v_i,i>n)=\]
\[P(n)=BP/(v_i,i<n)\]
\[B(n)=v_n^{-1}P(n)\]
\[k(n)=BP/(v_i,i\neq n)\]
\[K(n)=v_n^{-1}BP/(v_i,i\neq
n)\] \[E(n)=v_n^{-1}BP/(v_i,i> n).\] Note that inverting an
element in the coefficient ring of an $MU_{(p)}$-ring spectrum is
an instance of Bousfield localization in the homotopy category of
$MU_{(p)}$-modules. Further Bousfield localization preserves
algebra structures we conclude that all spectra listed above admit
structures of $MU_{(p)}$-algebras (and therefore $MU$-algebras).

Finally we need to mention the result of Hopkins and Miller, cf.
\cite{Rez}, which states that a suitable completion $E_n$ of the
spectrum $E(n)$ (which came to be popularly known as the Morava
$E$-theory) admits a unique structure of an $S$-algebra. A
generalization of this theorem is discussed in \cite{L2}.
Furthermore, in \cite{GH} Hopkins and Goerss proved using a more
advanced technology that $E_n$ in fact admits a unique structure
of a \emph{commutative} $S$-algebra.
\section{Spaces of multiplicative maps}
In this section we will consider the problem of computing the
homotopy type of the mapping space
$\Map_{\mathcal{C}^{ass}_R}(A,B)$ where $A$ is a cofibrant
$R$-algebra. Unexpectedly, it turns out that the space of based
loops on $\Map_{\mathcal{C}^{ass}_R}(A,B)$ is weakly equivalent to
an infinite loop space. More precisely, in \cite{L3} the following
theorem is proved.
\begin{theorem}\label{osn2}For a $q$-cofibrant algebra $A$ and a map of $R$-algebras
$f:A\rightarrow B$ the space $d$-fold  based loops
$\Omega^d(F_{R-alg}(A,B),f)$ is weakly equivalent to the space
$\Omega^\infty{\bf Der}_R(A,\Sigma^{-d}B)$.
\end{theorem}
The main ingredient in the proof of this theorem is the following
\begin{lem}\label{thy}
There is a weak equivalence of $R$-algebras $A^{S^d}$ and the
square-zero-extension $A\vee \Sigma^{-d}A$.\end{lem} Note that,
obviously, there is a weak equivalence of $R$-modules
\begin{equation}\label{der}\xymatrix{A^{S^d}\ar[r]&A\vee
\Sigma^{-d}A}\end{equation}
 The lemma is
proved by computing the topological derivations of $A^{S^d}$ and
showing that the map of (\ref{der}) could be improved to a
topological derivation of $A^{S^d}$ with values in $\Sigma^{-d}A$.

The proof of Theorem \ref{osn2} is actually very easy. We have the
following commutative diagram of spaces where both rows are
homotopy fibre sequences:
\[\xymatrix{{\mathcal Top}_\ast(S^d,\Map_{\mathcal{C}^{ass}_R}(A,B))
\ar[r]\ar[d]&{\mathcal
Top}(S^d,\Map_{\mathcal{C}^{ass}_R}(A,B))\ar[r]\ar[d]
&\Map_{\mathcal{C}^{ass}_R}(A,B)\ar@{=}[d]\\ ?\ar[r]&
\Map_{\mathcal{C}^{ass}_R}(A,B^{S^d})\ar[r]&\Map_{\mathcal{C}^{ass}_R}(A,B)}
\]
Here the horizontal rightmost arrows are both induced by the
inclusion of the base point into $S^d$. Since the right and the
middle vertical arrows are weak equivalences (even isomorphisms)
it follows that the map
 ${\mathcal Top}_\ast(S^d,\Map_{\mathcal {C}^{ass}_R}(A, B))\rightarrow ?$ is a weak equivalence. But
Lemma \ref{thy} tells us that the $R$-algebra  $B^{S^d}$ is weakly
equivalent as an $R$-algebra to $B\vee \Sigma^{-d}B$. In other
words the term $?$ is weakly equivalent to the topological space
of maps $A\rightarrow B\vee \Sigma^{-d}B$ which yield $f$ when
composed with the projection $B\vee \Sigma^{-d}B\lra B$. Therefore
$?$ is weakly equivalent to $\Omega^\infty{\bf
Der}_R(A,\Sigma^{-d}B)$ and our theorem is proved.
\begin{cor}\label{com}
For a $q$-cofibrant algebra $A$ and a map of $R$-algebras
$f:A\rightarrow B$ there is a bijection between sets
$\pi_d(\Map_{\mathcal{C}^{ass}_R}(A,B),f)$ and $Der^{-d}_R(A,B)$
for $d\geq 1$. If $d\geq 2$ then this bijection is an isomorphism
of abelian groups.
\end{cor}
One might wonder whether Theorem \ref{osn2} remains true in the
context of commutative $S$-algebras. The answer is no. The crucial
point is the weak equivalence of $S$-algebras $S\vee S^{-1}$ and
$S^{S^1}$.  It is clear that $\pi_0S\wedge_{S\vee S^{-1}}S $ is
the divided power ring. However N.Kuhn and M.Mandell proved that
$\pi_0S\wedge_{S^{S^1}}S $ is the ring of numeric polynomials.
Therefore  $S\vee S^{-1}$ and $S^{S^1}$ cannot be weakly
equivalent as commutative $S$-algebras.

Another point of view is afforded by Mandell's theorem \cite{Man}
which states that the homotopy category of connected $p$-complete
nilpotent spaces of finite $p$-type can be embedded into the
category of $E_\infty$ algebras over $\bar{\mathbb{F}}_p$ as a
full subcategory. Here $\bar{\mathbb F}_p$ is the algebraic
closure of the field $\mathbb{F}_p$. This embedding is via the
singular cochains functor $C^*(-,\mathbb{F}_p)$. The category of
$E_\infty$ algebras over $\bar{\mathbb{F}}_p$ is Quillen
equivalent to the category of commutative
$H\bar{\mathbb{F}}_p$-algebras. Consider an arbitrary space $X$
(connected, $p$-complete, nilpotent and of finite $p$-type) and
the corresponding $H\bar{\mathbb{F}}_p$-algebra $\tilde X$ which
we assume to be cofibrant. Then the space of maps in
$\mathcal{C}^{comm}_{H\bar{\mathbb{F}}_p}$ from $\tilde X$ into
$H\bar{\mathbb{F}}_p$ will be weakly equivalent to $X$. This shows
that the space of commutative $R$-algebra maps could have an
essentially arbitrary homotopy type.

However, there is a context in which the analogue of Theorem
\ref{osn2} holds for commutative $R$-algebras. Let us suppose that
our ground $S$-algebra $R$ is rational, that is, local with
respect to $H\mathbb{Q}$. Let $A$ and $B$ be commutative
$R$-algebras where $A$ is cofibrant as before. Then it is easy to
see that $A^{S^1}$ and $A\vee \Sigma^{-1}A$ are weakly equivalent
as \emph{commutative} $R$-algebras. Indeed, $A\vee \Sigma^{-1}A$
is weakly equivalent as a commutative $R$-algebra $A\wedge
R[S^{-1}]$ where $R[S^{-1}]$ is the free commutative $R$-algebra
on the $R$-module $S^{-1}$. We have
$\pi_*A^{S^1}=\Lambda_{A_*}(x)$ where the exterior generator $x$
has degree $-1$. The element $x$ determines a map of $R$-modules
${S^1}\lra A^{S^1}$ and therefore a map of $R$-algebras $A\vee
\Sigma^{-1}A\simeq A\wedge R[S^{-1}]\lra A^{S^1}$ which is a weak
equivalence. Further there are the following equivalences of
commutative $R$-algebras:
\[A^{S^d}\cong A^{S^1\wedge S^{d-1}}\cong
(A^{S^1})^{S^{d-1}}\simeq (A\vee \Sigma^{-1}A)^{S^{d-1}}\simeq
A\vee \Sigma^{-d}A.\] So we in this case we have the exact
analogue of Lemma \ref{thy}. The remainder of the proof is the
same and we obtain
\begin{theorem}\label{osn3}Let $R$ be a rational commutative $S$-algebra,
Then for a $q$-cofibrant $R$-algebra $A$ and a map of $R$-algebras
$f:A\rightarrow B$ the space $d$-fold based loops
$\Omega^d(\Map_{\mathcal{C}^{comm}_R},f)$ is weakly equivalent to
the space $\Omega^\infty{\bf TAQ}_R(A,\Sigma^{-d}B)$. Furthermore
there is a bijection between sets
$\pi_d(\Map_{\mathcal{C}^{comm}_R}(A,B),f)$ and $TAQ^{-d}_R(A,B)$
for $d\geq 1$. If $d\geq 2$ then this bijection is an isomorphism
of abelian groups.
\end{theorem}
For $R=H\mathbb{Q}$, the Eilenberg-MacLane spectrum the category
of commutative $R$-algebras is Quillen equivalent to the category
of commutative differential graded algebras over $\mathbb{Q}$. Via
this equivalence Theorem \ref{osn3} translates into a statement in
rational homotopy theory which does not seem to have appeared in
the literature.
\section{Moduli spaces}
For an $R$-module $X$ it is natural to look at the set of all
nonisomorphic $R$-algebra structures on $X$. (We assume that $X$
supports at least one $R$-algebra structure.) This set is actually
$\pi_0$ of the Dwyer-Kan \emph{classification space}, cf.
\cite{DK} which we will now describe.

Consider the subcategory $W$ of $\mathcal{C}^{ass}_R$ consisting
of those $R$-algebras which are weakly equivalent to $X$ as
$R$-modules. The morphisms of $W$ are the weak equivalences of
$R$-algebras. Then the \emph{moduli space} $\mathcal{M}(X)$ is by
definition the nerve of $W$.

For a cofibrant $R$-algebra $A$ we will denote by $haut(A)$ the
topological monoid of homotopy autoequivalences of $A$. For a
noncofibrant $A$ we take for $haut(A)$ the monoid of homotopy
autoequivalences of its cofibrant approximation. Thus $haut(A)$ is
well defined up to homotopy. We can form its classifying space
$Bhaut(A)$. Then the results of Dwyer and Kan imply the following
\begin{theorem}\label{dwyerkan} There is a weak equivalence of spaces
\[\mathcal{M}(X)\simeq \coprod Bhaut(A)\]
where the disjoint union is taken over the set of connected
components of $W$.
\end{theorem}
We assumed from the beginning the existence of at least one
$R$-algebra structure on $X$, that gives the space
$\mathcal{M}(X)$ a base point. Let us denote the corresponding
$R$-algebra by $A$. The connected component of $A$ in
$\mathcal{M}(X)$ is, therefore, weakly equivalent to $Bhaut(A)$.

Note that it is possible to consider other moduli spaces for $X$.
For example, we could insist that the homotopy $X_*$ has a fixed
ring structure, or that $X$ itself has a fixed $R$-ring spectrum
structure up to homotopy. Each of these choices affects $\pi_0$ of
our moduli space but not the connected components themselves, and
therefore, not the higher homotopy groups.

Assume, without loss of generality, that $A$ is a cofibrant
$R$-algebra. Clearly $haut(A)$ is a disjoint union summand of the
space $\Map_{\mathcal{C}^{ass}_R}(A,A)$ of all multiplicative
self-maps of $A$. Taking the identity map $A\lra A$ as the base
point in both $haut(A)$ and $\Map_{\mathcal{C}^{ass}_R}(A,A)$ we
see that
\[\Omega haut(A)\simeq \Omega\Map_{\mathcal{C}^{ass}_R}(A,A).\]
But we know from the previous section that there is a weak
equivalence of spaces
\[\Omega\Map_{\mathcal{C}^{ass}_R}(A,A)\simeq
\Omega\Omega^{\infty}{\bf Der}_R(A,A).\] This gives the following
\begin{theorem}\label{moduliderivation}The space $\Omega^2\mathcal{M}(X)$ of two-fold
loops on the moduli space $\mathcal{M}(X)$ is weakly equivalent to
$\Omega\Omega^{\infty}{\bf Der}_R(A,A)$. In particular it is an
infinite loop space.
\end{theorem}
In particular, we have an isomorphism
\[\pi_i\mathcal{M}(X)\cong Der_R^{1-i}(A,A)\]
for $i\geq 2$. The computation of $\pi_1$ and $\pi_0$ is, of
course, a completely different story.

What about the moduli space of \emph{commutative} $R$-algebra
structures? Theorem \ref{dwyerkan} has an obvious analogue in the
context of commutative $R$-algebras. However the analogue of
Theorem \ref{moduliderivation} falls through because for a
commutative $R$-algebra $A$ the space of self-maps of $A$ in
$\mathcal{C}^{comm}_{R}$ is not in any obvious way related to
${\bf TAQ}_R(A,A)$. However Hopkins and Goerss proved in \cite{GH}
that $E_n$, the Morava $E$-theory spectrum mentioned at the end of
Section 7 admits a unique structure of a commutative $S$-algebra
and the space of its self-maps is homotopically discrete with the
set of connected components being equal to $S_n$, the Morava
stabilizer group. In our interpretation it means that the
appropriate moduli space of $E_n$ is weakly equivalent to $BS_n$,
the classifying space of $S_n$.


\begin{thebibliography}{40}
\bibitem{Bas} M. Basterra, {\em Andr\'e-Quillen cohomology of commutative $S$-algebras}.
J. Pure Appl. Algebra 144 (1999), no. 2, 111--143.
\bibitem{BM}M. Basterra \& R. McCarthy, {\em $\Gamma$-homology, topological Andr\'e-Quillen
homology and stabilization.} Topology Appl. 121 (2002),
no.3,551--566.
\bibitem{BaL} A.Baker \& A.Lazarev, {\em Topological Hochschild cohomology and generalized
Morita equivalence}, preprint, math.AT/0209003.
\bibitem{Bok} M.B\"okstedt. {\em Topological Hochschild homology of
$\mathbb{Z}$ and $\mathbb{Z}_p$}, preprint.
\bibitem{QC} J. Cuntz \& D. Quillen, {\em Algebra extensions and nonsingularity}.
J. Amer. Math. Soc. 8 (1995), no. 2, 251--289.
\bibitem{DK} W. Dwyer \& D. Kan,{\em A classification theorem for diagrams of
simplicial sets.} Topology 23 (1984), no. 2, 139--155.
\bibitem{EKMM} A.
Elmendorf, I. Kriz, M. Mandell \& J. P. May, {\em Rings, modules,
and algebras in stable homotopy theory}. Mathematical Surveys and
Monographs {\bf 47} (1997).
\bibitem{Fra} V. Franjou \&  J. Lannes, L. Schwartz, {\em Autour de
la cohomologie de Mac Lane des corps finis.}   Invent. Math. 115
(1994), no. 3, 513--538.
\bibitem{Goe} P.Goerss, {\em Associative $MU$-algebras}, preprint 2001.
\bibitem{GH} P.Goerss \& M.Hopkins {\em Realizing commutative ring spectra as
$E_\infty$ ring spectra}, preprint.
\bibitem{Goo1}T. Goodwillie, {\em Calculus I. The first derivative of pseudoisotopy theory.}
$K$-Theory 4 (1990), no. 1, 1--27.
\bibitem{Goo3} T. Goodwillie, {\em Calculus 3. The Taylor series of
homotopy functors}, preprint.
\bibitem{Hov} M. Hovey, {\em Model categories.}
Mathematical Surveys and Monographs, 63. American Mathematical
Society, Providence, RI, 1999.
\bibitem{HSS}M. Hovey, B. Shipley \& J. Smith, {\em Symmetric spectra.}
J. Amer. Math. Soc. 13 (2000), no. 1, 149--208.
\bibitem{Kriz} I.Kriz. {\em Towers of $E_\infty$-ring spectra with applications
to $BP$}, preprint, 1993.
\bibitem{L1} A. Lazarev, {\em Homotopy theory of $A\sb \infty$ ring spectra and
applications to $M{\rm U}$-modules}. $K$-Theory 24 (2001), no. 3,
243--281.
\bibitem{L2} A.Lazarev. {\em Towers of $MU$-algebras and a generalized
Hopkins-Miller theorem}, to appear in the Proc. of the London
Math.Soc.
\bibitem{L3} A.Lazarev, {\em Spaces of multiplicative maps between
highly structured ring spectra}, Proceedings of the Isle of Sky
Conference on Algebraic Topology, to appear.
\bibitem{Mad}A. Lindenstrauss \& I. Madsen, {\em Topological Hochschild homology of
number rings.} Trans. Amer. Math. Soc. 352 (2000), no. 5,
2179--2204.
\bibitem{Madsen} I. Madsen, {\em Algebraic $K$-theory and traces}.
Current developments in mathematics, 1995 (Cambridge, MA),
191--321, Internat. Press, Cambridge, MA, 1994.
\bibitem{Man} M. Mandell, {\em $E\sb \infty$ algebras and $p$-adic homotopy theory}.
Topology 40 (2001), no. 1, 43--94.
\bibitem{M1} J.P May,  {\em $E\sb{\infty }$ ring spaces
and $E\sb{\infty }$ ring spectra}. With contributions by F. Quinn,
N. Ray, and J. Tornehave. Lecture Notes in Mathematics, Vol. 577.
Springer-Verlag, Berlin-New York, 1977.
\bibitem{May} J.P. May, {\em The geometry of iterated loop spaces}.
Lectures Notes in Mathematics, Vol. 271. Springer-Verlag,
Berlin-New York, 1972.
\bibitem{MS} J. McClure \& R. Staffeldt, {\em On the topological Hochschild
homology of $b{\rm u}$}. I. Amer. J. Math. 115 (1993), no. 1,
1-45.
\bibitem{PR} T.Pirashvili \& B. Richter, {\em Robinson-Whitehouse complex and
stable homotopy.} Topology 39 (2000), no. 3, 525--530.
\bibitem{Q} D.Quillen, {\em Homotopical algebra}.
Lecture Notes in Mathematics, No. 43. Springer-Verlag, Berlin-New
York 1967.
\bibitem{Rez} C. Rezk, Notes on the Hopkins-Miller Theorem.
{\em Homotopy theory via algebraic geometry and group
representations} (Evanston, IL, 1997), 313--366, Contemp. Math.,
220, Amer. Math. Soc., Providence, RI, 1998.
\bibitem{Ric} B. Richter, {\em An Atiyah-Hirzebruch spectral sequence
for topological Andr\'e-Quillen homology.} J. Pure Appl. Algebra
171 (2002), no. 1, 59--66.
\bibitem{RW} A. Robinson \& S. Whitehouse, {\em Operads and $\Gamma$-homology
of commutative rings.} Math. Proc. Cambridge Philos. Soc. 132
(2002), no. 2, 197--234.
\bibitem{Rob1}A.Robinson, {\em Obstruction theory and the strict associativity
of Morava $K$-theories.} Advances in homotopy theory 143--152,
London Math. Soc. Lecture Note Ser., 139, Cambridge Univ. Press,
Cambridge, 1989.
\bibitem{Rob2} A.Robinson, {\em Gamma homology, Lie representations and
$E_\infty$ multiplications}, preprint available at
\texttt{http://www.maths.warwick.ac.uk/$\sim$car/}
\bibitem{Str} N.P.Strickland, {\em Products on $MU$-modules.}
Trans. Amer. Math. Soc. 351 (1999), no. 7, 2569--2606.
\end{thebibliography}
\end{document}